\newif\ifsmfart
\IfFileExists{smfart.cls}
   {\documentclass[12pt,english]{smfart}
    \IfFileExists{smfenum.sty}{\usepackage{smfenum}}{}
    \usepackage{bull}
    \smfarttrue}
   {\message{^^J*** It would be better to typeset this file with smfart.cls ***^^J^^J}
    \documentclass[12pt,a4paper]{amsart}}

\usepackage{amssymb}

 \let\mathsec\mathsf
 \IfFileExists{mathrsfs.sty}
  {\usepackage{mathrsfs}\let\mathcal\mathscr}
  {\let\mathscr\mathcal}

\advance\headheight 2pt
\calclayout
\allowdisplaybreaks[3]

\numberwithin{equation}{section}
\makeatletter

\let\c@equation\c@subsection
\let\cl@equation\cl@subsection
\def\theequation{\thesection.\arabic{equation}}

\def\l@table{\@tocline{0}{3pt plus2pt}{0pt}{}{\itshape}}

\makeatother

\theoremstyle{plain}
\newtheorem{prop}[subsection]{Proposition}

\newtheorem{theo}[subsection]{Theorem}
\newtheorem{coro}[subsection]{Corollary}

\newtheorem{lemm}[subsection]{Lemma}

\newtheorem{defi}[subsection]{Definition}

\theoremstyle{definition}

\theoremstyle{remark}
\newtheorem{rema}[subsection]{Remark}

\newtheorem{exam}[subsection]{Example}

\def\A{{\mathbf A}}
\def\Afin{{{\mathbf A}_{\textrm{fin}}}}
\def\C{{\mathbf C}}

\def\G{{\mathbf G}}
\def\ga{{\mathbf G}_a}
\def\gm{{\mathbf G}_m}
\def\N{{\mathbf N}}
\def\P{{\mathbf P}}
\def\Q{{\mathbf Q}}
\def\R{{\mathbf R}}

\def\Z{{\mathbf Z}}

\def\Res{{\operatorname{Res}}}
\def\res{{\operatorname{res}}}

\def\dx{{\mathrm d}{\mathbf x}}

\let\ra\rightarrow

\let\epsilon\varepsilon \let\eps\epsilon
\let\epsilon\varepsilon
\let\phi\varphi
\let\emptyset\varnothing
\let\leq\leqslant
\let\geq\geqslant
\let\le\leq
\let\ge\geq

\def\eff{{\text{\upshape eff}}}
\def\abs#1{\left\lvert{#1}\right\rvert}
\def\norm#1{\left\|{#1}\right\|}

\def\DeclareMathOperator#1#2{\def #1{\operatorname{#2}}}
\DeclareMathOperator{\red}{red} 
\DeclareMathOperator{\Re}{Re} 
\DeclareMathOperator{\Im}{Im} 
\DeclareMathOperator{\Pic}{Pic}
\DeclareMathOperator{\Gal}{Gal}

\DeclareMathOperator{\Spec}{Spec}

\DeclareMathOperator{\rk}{rk}

\DeclareMathOperator{\div}{div}

\DeclareMathOperator{\vol}{vol}

\def\dx{\,{\mathrm d\mathbf x}}

\title[Points of bounded height on compactifications of vector groups]
      {On the distribution \\ of points of bounded height \\
       on equivariant compactifications \\ of vector groups }

\author{Antoine Chambert-Loir}
\address{Institut de math{\'e}matiques de Jussieu, Boite 247 \\ 
4 place Jussieu \\ 
F-75252 Paris Cedex 05 }
\email{chambert@math.jussieu.fr}

\author{Yuri Tschinkel}
\address{Department of Mathematics, Princeton University\\
Princeton, (NJ) 08544-1000,  U.S.A. }
\email{ytschink@math.princeton.edu}

\expandafter\ifx\csname pdfoutput\endcsname\relax
\else
 \usepackage[
             pdftitle={Points of bounded heights...},
             pdfauthor={A. Chambert-Loir \& Y. Tschinkel}]
             {hyperref}
\fi

\begin{document}
\def\smfandname{\&}
\date{First version submitted on the arXiv, May 2, 2000.
Second version, \today}

\begin{abstract}
We prove asymptotic formulas  for the number
of rational points of bounded height
on smooth equivariant compactifications of the affine space.
\end{abstract}
 
\ifsmfart 
 \begin{altabstract}
 Nous {\'e}tablissons un d{\'e}veloppement asymptotique
 du nombre de points rationnels de hauteur born{\'e}e sur
 les com\-pac\-ti\-fi\-ca\-tions {\'e}qui\-va\-riantes lisses
 de l'espace affine.
 \end{altabstract}
\fi

\maketitle


\setcounter{section}{-1}
\section{Introduction}

A theorem of Northcott asserts that for any real number $B$
there are only finitely many rational points in 
the projective space $\P^n$ with {\em height}
smaller than~$B$. An asymptotic formula for
this number (as $B$ tends to infinity) has been proved
by Schanuel~\cite{schanuel79}.
Naturally, it is interesting to consider 
more general projective varieties
and there are indeed a number of results in this direction.
The techniques employed can be grouped in three main classes:
\begin{itemize}
\item the classical circle method in analytic number theory
permits to treat complete intersections of small degree in
projective spaces of large dimension
(cf.~for example~\cite{birch62});
\item harmonic analysis on adelic points of reductive
groups leads to results for toric varieties~\cite{batyrev-t98b},
flag varieties~\cite{franke-m-t89}
and horospherical varieties~\cite{strauch-t99};
\item elementary (but nontrivial) methods for
del Pezzo surfaces of degree 4 or 5
cf.\ de~la~Bret{\`e}che~\cite{breteche2000},
and for cubic surfaces
(Salberger, Swinnerton-Dyer, Heath-Brown, cf.~\cite{proc-peyre98}
and the references therein;
\end{itemize}

\medskip

This research has been stimulated by a conjecture put
forward by Batyrev and Manin.
They proposed
in~\cite{batyrev-m90}
an interpretation of the growth
rate in terms of the mutual positions of the class of the line 
bundle giving the projective embedding,
the anticanonical class and the cone of effective divisors in the
Picard group of the variety.  
Peyre refined this conjecture in~\cite{peyre95}
by introducing an adelic Tamagawa-type number
which appears as the leading constant in 
the expected asymptotic formula for the anticanonical embedding.
Batyrev and the second author proposed an interpretation
of the leading constant for arbitrary ample line bundles,
see~\cite{batyrev-t98}. 

\medskip

In this paper we consider
a new class of varieties, namely
\emph{equivariant compactifications of vector groups}.
On the one hand, we can make use of 
harmonic analysis on the adelic points of the group. 
On the other hand, such varieties
have a rich geometry. In particular, in contrast
to flag varieties and toric varieties, they admit geometric deformations,
and in contrast to the complete intersections treated by the
circle method their Picard group can have arbitrarily high rank.
Their geometric classification is a 
difficult open problem already in dimension
3 (see \cite{hassett-t99}).  

A basic example of such algebraic varieties is of course
$\P^n$ endowed with the action of $\ga^n\subset\P^n$ by translations.
A class of examples is provided by the following
geometric construction: 
Take $X_0=\P^n$ endowed with
the translation action of $\ga^n$. Let $Y$ be a smooth subscheme
of $X_0$ which is contained in the hyperplane at infinity.
Then, the blow-up $X=\operatorname{Bl}_Y(X_0)$ contains the isomorphic
preimage of $\ga^n$ and the action of $\ga^n$ lifts to~$X$.
The rank of $\Pic(X)$ is equal to the number of irreducible
components of~$Y$. Using equivariant resolutions of singularities,
one can produce even more complicated examples, although less explicitly.

Our first steps towards this paper
are detailed in~\cite{chambert-loir-t2000}
and~\cite{chambert-loir-t2000b}. There we studied 
the cases $n=2$ with $Y$ a finite union of $\mathbf Q$-rational points
and $n>2$ with $Y$ a smooth hypersurface contained in the
hyperplane at infinity.

\bigskip

We now describe the main theorem of this article.
\emph{Let $X$ be an equivariant
compactification of the additive group $\ga^n$ over
a number field~$F$.}
Unless explicitly stated, we shall always \emph{assume that~$X$
is smooth and projective.} The boundary divisor
$D=X\setminus\ga^n$ is a sum of irreducible components
$D_\alpha$ ($\alpha\in\mathcal A$). We do not assume that
they are geometrically irreducible.


The Picard group of $X$ is free and has a canonical
basis given by the classes of $D_{\alpha}$. The cone 
of effective divisors 
consists of  the divisors $\sum_{\alpha \in {\mathcal A}} 
d_{\alpha}D_{\alpha}$
with $d_{\alpha}\geq 0$ for all $\alpha$.  
Denote by  $K_X^{-1}$ the anticanonical line bundle on $X$ and by
$\rho=(\rho_\alpha)$ its class in $\Pic(X)$.

Let $\lambda=(\lambda_{\alpha})$ 
be a class contained in the interior of the 
cone of effective divisors $\Lambda_{\eff}(X)\subset \Pic(X)$ 
and $\mathcal L_{\lambda}$ 
the corresponding line bundle, equipped 
with a smooth adelic metric (see~\ref{defi.adelic}
for the definition). With the above notations, we have $\lambda_\alpha>0$
for all~$\alpha$.
Denote by $H_{\mathcal L_{\lambda}}$ the associated
exponential height on $X(F)$. 
Let 
$a_{\lambda}= \max(\rho_{\alpha}/\lambda_{\alpha})$
and let $b_{\lambda}$ be the cardinality of 
\[
{\mathcal B}_{\lambda}=\{ \alpha\in {\mathcal A}\, ;
\, \rho_{\alpha}=a_{\lambda}\lambda_{\alpha}\}.
\]
Put 
\[
c_{\lambda} = 
\prod_{\alpha\in {\mathcal B}_{\lambda}}\lambda_\alpha^{-1}.
\]
For example, one has $a_{\rho}=1$ and $b_{\rho}=\rk\Pic(X)$.
We denote by  $\tau({\mathcal K}_X)$ 
the Tamagawa number
as defined by Peyre in~\cite{peyre95}. 

\begin{theo}\label{theo.main}
\textup{a)}
The series
\[ Z_{\lambda}(s) = \sum_{x\in \ga^n(F)} H_{\mathcal L_{\lambda}}(x)^{-s} \]
converges absolutely and uniformly for $\Re(s)>a_{\lambda}$ 
and has a meromorphic
continuation to $\Re(s)>a_{\lambda}-\delta$ for some $\delta >0$, 
with a unique pole at $s=a_{\lambda}$
of order $b_{\lambda}$. Moreover, $Z_{\lambda}$ has polynomial growth in
vertical strips in this domain.

\textup{b)}
There exist positive
real numbers $\tau_{\lambda},\delta'$ 
and a polynomial $P_{\lambda}\in\R[x]$ of degree 
$b_{\lambda}-1$ with leading coefficient 
\[
c_{\lambda}\tau_{\lambda}/(b_{\lambda}-1)!
\]
such that the number $N({\mathcal L}_{\lambda},B)$ of $F$-rational 
points in $\ga^n$ with height $H_{\mathcal L_{\lambda}}$ 
smaller than $B$ satisfies
\[ N({\mathcal L}_{\lambda}, B)
= B^{a_{\lambda}} P_{\lambda}(\log B) + O(B^{a_{\lambda}-\delta'}), \]
for $B\ra \infty$.

\textup{c)}
For $\lambda =\rho$ we have 
$\tau_{\rho}=\tau({\mathcal K}_X)$. 
\end{theo}


\begin{rema}
Granted the smoothness assumption on the adelic metric,
note that the normalization of the height in~Theorem~\ref{theo.main}
can be arbitrary.
By a theorem of Peyre~\cite[\S\,5]{peyre95}, this means that
the points of bounded anticanonical height are equidistributed 
with respect to the Tamagawa measure
(compatible with the choice of the height function)
in the adelic space $X(\A_F)$.
Let us explain this briefly.
Fix a smooth adelic metric on the anticanonical line bundle,
this defines a height function~$H$. Let $\mathrm d\tau_H$
be the renormalized Tamagawa measure on $X(\A_F)$.
A smooth positive function $f$ on $X(\A_F)$
determines another height function, namely $H'=fH$.
Applied to such~$H'$, the main theorem implies that 
\[  \lim_{B\ra+\infty} 
     \frac{(r-1)!}{c_{\rho}} \frac{1}{B (\log B)^{r-1}}
     \sum_{\substack{x\in \ga^n(F) \\ H(x)\leq B}} f(x)
    = \int_{X(\A_F)} f(\mathbf x) \mathrm d\tau_H(\mathbf x).
\] 
\end{rema}

\begin{rema}
Theorem~\ref{theo.main}
implies that the open subset $\ga^n$ does not contain any
accumulating subvarieties. 
\end{rema}

The proof proceeds as follows.

First, we extend the height function to the
adelic space $\ga^n(\A_F)$. Next we apply the additive
Poisson formula and find a representation of $Z_{\lambda}(s)$ as
a sum over the characters of $\ga^n(\A_F)/\ga^n(F)$
of the Fourier transforms of $H$.

The Fourier transforms of~$H$ decompose as products over all
places $v$ of ``global'' integrals on $X(F_v)$ which are
reminiscent of Igusa zeta functions.

The most technical part of the paper is devoted to the evaluation
of these products: meromorphic continuation and control of
their growth in vertical strips.
For this, we need to consider the special case
where $D$ has strict normal crossings:
this means that over an algebraic closure of~$F$,
$D$ is a sum of smooth irreducible divisors meetings transversally.
Then, at almost all nonarchimedean places,
we compute explicitly the local Fourier transforms in terms of
the reduction of $X$ modulo the corresponding prime.  
The obtained formulas resemble Denef's
formula in~\cite{denef87} for Igusa's local zeta function.
For the remaining nonarchimedean places we find estimates.
This leads to a proof of the meromorphic continuation of
the Fourier transforms at each character.

To consider the general cas, where $D$ is not assumed to have
strict normal crossings,
we then introduce a proper modification $\pi\colon\tilde X\ra X$
which is a composition of equivariant blow-ups with smooth centers
lying on~$D$ such that 
$\tilde X$ is a smooth projective equivariant compactification
of~$\ga^n$ whose boundary divisor $\tilde D$ has strict normal crossings:
over an algebraic closure of~$F$, $\tilde D$ is a sum of smooth
irreducible divisors meeting transversally.
Considering $\pi^*\lambda$, all the computations shall be made on $\tilde X$.
We also prove in Lemma~\ref{lemma.singularities}
that $a_{\pi^*\lambda}=a_{\lambda}$
and $b_{\pi^*\lambda}=b_{\lambda}$.

In the Poisson formula,
invariance properties of the height   
reduce
the summation over $\ga^n(F)$ to one over a lattice.
The meromorphic continuation in part a) of the theorem
follows then from additional estimates for the Fourier
transforms at the infinite places. 

At this stage one has a meromorphic continuation of $Z_{\lambda}(s)$
to the domain $\Re(s)>a_{\lambda}-\delta$ for 
some $\delta>0$, with a single pole at $s=a_{\lambda}$
whose order is less or equal than $b_{\lambda}$.
It remains  to check that the order of the pole is exactly $b_{\lambda}$;
we need to prove that the limit
\[
\lim_{s\ra a_{\lambda}} Z_{\lambda}(s)(s-a_{\lambda})^{b_{\lambda}}
\]
is strictly positive.

For $\lambda=\rho$ this is more or less straightforward: 
the main term is given by the summand corresponding
to the trivial character and
the Tamagawa number defined by Peyre appears naturally in the limit.
For other $\lambda$ we use the Poisson formula
again and relate the limit to an integral of the height over
some subspace which is shown to be strictly positive.

Part b) follows  by a Tauberian theorem.

%
%
%
%

\subsection*{Acknowledgments} 
The work of the second author was partially supported by
the NSA.


\section{Geometry} \label{section.geometry}

Let $X$ be a smooth projective equivariant compactification
of the additive group $\ga^n$ 
of dimension $n$ over a field $F$: 
$X$ is a (smooth, projective) algebraic variety
endowed with an action of~$\ga^n$
with a dense orbit isomorphic to~$\ga^n$.
We will always assume that  
$X$ is smooth and projective. Let $D=X\setminus\ga^n$
denote the boundary divisor. We have a decomposition
in irreducible components
\[  D = \bigcup_{\alpha\in{\mathcal A}} D_\alpha \]
where $\mathcal A$ is a 
finite set and for each $\alpha\in \mathcal A$,
$D_\alpha$ is an integral divisor in $X$. 
We do not assume that the $D_\alpha$ are geometrically irreducible,
nor that they are smooth.

There is however a similar description over any extension of~$F$.
In particular, let $\bar F$ be a separable closure of~$F$
and let
\[ D_{\bar F} = \bigcup_{\alpha\in {\mathcal A}_{\bar F}}  D_{\bar F,\alpha} \]
be the decomposition of $D_{\bar F}$ in irreducible components.
The natural action of $\Gamma_F=\Gal(\bar F/F)$ on $\bar X$
induces an action of $\Gamma_F$ on $\mathcal A_{\bar F}$ such
that for any $\alpha\subset\mathcal A_{\bar F}$, $g(D_\alpha)=D_{g(\alpha)}$.

\begin{prop} \label{prop.picard}
One has natural isomorphisms of $\Gamma_F$-modules
(resp.~$\Gamma_F$-monoids)
\[ \bigoplus_{\alpha\in \mathcal A} \Z D_\alpha \ra \Pic(\bar X)
\quad\text{and}\quad
\bigoplus_{\alpha\in \mathcal A} \N D_\alpha \ra \Lambda_{\eff}(\bar X), 
\]
where  $\Lambda_{\eff}(X)$ is the monoid
of classes of effective divisors of $X$. 
\end{prop}

Subsequently, we identify the divisors 
$D_{\alpha}$ (and the corresponding line bundles) 
with their classes in the Picard group.

\begin{proof}
These maps are equivariant under the action of $\Gamma_F$. Hence, it
remains to show injectivity and surjectivity.

For this we may assume $\bar F=F$.
Let $\mathcal L$ be a line bundle on $X$.
As $X$ is smooth and $\Pic(\ga^n)=0$, there is a divisor $D$ in $X$ not
meeting $\ga^n$ such that $\mathcal L\simeq\mathcal O_X(D)$.
Such a divisor $D$ is a sum $\sum_{\alpha\in A} n_\alpha D_\alpha$.
Moreover, such a $D$ is necessarily unique.
If $\mathcal O_X(D)\simeq \mathcal O_X(D')$ for
$D=\sum_\alpha n_\alpha D_\alpha$ and
$D'=\sum_\alpha n'_\alpha D_\alpha$, then the canonical
rational section $\mathsec s_D/\mathsec s_{D'}$  of $\mathcal
O_X(D-D')=\mathcal O_X$ is a rational function on $X$ without
zeroes nor poles on $\ga^n$. 
Rosenlicht's lemma~(see lemma~\ref{lemm.rosenlicht} below) implies
that $\mathsec s_D/\mathsec s_{D'}$ is constant, so that $D=D'$.

We remark that any effective cycle $Z$ on $X$
is rationally equivalent to
a cycle that does not meet $\ga^n$. Indeed, if $t$ is the parameter
of a subgroup of $\ga^n$ isomorphic to $\ga$, we can consider
the specialization of the cycles $t+Z$ when $t\ra\infty$.
\end{proof}

\begin{lemm}[Rosenlicht] \label{lemm.rosenlicht}
If $f\in F(X)$ has neither zeroes nor poles on $\ga^n$,
then $f\in F^*$.
\end{lemm}

\begin{coro}
{\upshape 1)} The $\Gamma_F$-module $\Pic(\bar X)$ is a permutation module.
In particular, \hskip0pt plus3cm\penalty-30\hskip0pt plus-3cm
$\mathrm H^1(\Gamma_F,\Pic(\bar X))=0$.

{\upshape 2)} $\Pic(X)=\Pic(\bar X)^{\Gamma_F}$ is a free $\Z$-module
of finite rank equal to the number of
$\Gamma_F$-orbits in~$\mathcal A_{\bar F}$.
\end{coro}

For any $\alpha\in\mathcal A$, we shall denote by
by $F_\alpha$ the algebraic closure of $F$ in the function
field of $D_\alpha$.
Over $F_\alpha$, the irreducible components
of $D_\alpha$ are  geometrically irreducible.
We also denote by $\zeta_{F_\alpha}$ the Dedekind zeta function
of the number field~$F_\alpha$.

\bigskip

Let $f$ be a nonzero linear form on $\ga^n$ viewed as an element of $F(X)$.
Its divisor can be written as
\[ \div(f) = E(f) - \sum_{\alpha\in \mathcal A} d_\alpha(f) D_\alpha \]
where $E(f)$ is the unique irreducible component of $\{f=0\}$ that
meets $\ga^n$ and $d_\alpha(f)$ are integers.
(The divisor $E(f)$ can also be seen 
as the closure in $X$ of the hypersurface
of $\ga^n$ defined by $f$.)

Since $E(f)$ is rationally equivalent to 
$\sum_{\alpha\in A} d_\alpha(f) D_\alpha$, 
the preceding proposition implies the following lemma:
\begin{lemm} 
For any nonzero linear form $f$ and any $\alpha\in\mathcal A$,
one has $d_\alpha(f)\geq 0$.
\end{lemm}

If $\mathbf a\in\ga^n(F)$, let $f_{\mathbf a}=\langle \cdot,\mathbf a\rangle$
be the associated linear form and define
\begin{equation}
\mathcal A_0(\mathbf a)= \{\alpha\,;\,d_{\alpha}(f)=0\}\quad\text{and}\quad
\mathcal A_0(\mathbf a)= \{\alpha\,;\,d_{\alpha}(f)=1\}.
\end{equation}

\begin{prop} \label{prop.invsection}
Let $X$ be a normal equivariant compactification of $\ga^n$ over $F$. 
Every line bundle $\mathcal L$ on $X$ admits a
unique $\ga^n$-linearization.
If $\mathcal L$ is effective then $\mathrm H^0(X,\mathcal L)$ has a unique
line of $\ga^n$-invariant sections.
\end{prop}
\begin{proof}
Since $\ga^n$ has no nontrivial characters,
it follows from Proposition 1.4 and from the proof of Proposition 
1.5 in~\cite{mumford-f-k94}, Chapter~1,
that any line bundle on $X$ admits a unique $\ga^n$-linearization.

Assume that $\mathrm H^0(X,\mathcal L)\neq 0$ and 
consider the induced action of $\ga^n$ on
the projectivization $\P(\mathrm H^0(X,\mathcal L))$.
Borel's fixed point theorem (\cite{borel91}, Theorem 10.4)
implies that there
exists a nonzero section $\mathsec s \in \mathrm H^0(X,\mathcal L)$
such that the line $F\mathsec s$ is fixed under this action.
As $\ga^n$ has no nontrivial characters, $\mathsec s$ itself is
fixed. The divisor $\div(\mathsec s)$ is
$\ga^n$-invariant; therefore, $\div(\mathsec s)$
does not meet $\ga^n$ and is necessarily a sum $\sum d_\alpha D_\alpha$
such that $\mathcal L\simeq \mathcal O_X(\sum d_\alpha D_\alpha)$.
Because of Proposition~\ref{prop.picard}, every other such 
section will be proportional to $\mathsec s$.
\end{proof}

\begin{rema}
We observe that if $D=\sum d_\alpha D_\alpha$
for integers $d_\alpha\geq 0$,
then the canonical section $\mathsec s_D$ of $\mathcal O_X(D)$ is
$\ga^n$-invariant.
\end{rema}

\section{Vector fields}\label{section.vector}

We now recall some facts concerning vector fields on equivariant
compactifications of algebraic groups.
Let $G$ be a connected algebraic group over $F$ and $\mathfrak g$ its Lie
algebra of invariant vector fields. 
Let $X$ be a smooth equivariant compactification of $G$.
Denote by $D=X\setminus G$ the boundary.
We assume that $D$ is a divisor with
strict normal crossings.
Let $\mathcal T_X$ be the tangent bundle of $X$.
Evaluating a vector field at the neutral 
element $\mathbf 1$ of $G$ induces a
``restriction map''
\[ \mathrm H^0(X,\mathcal T_X) 
\rightarrow \mathcal T_{X,\mathbf 1}=\mathfrak g. \]
Conversely, given $\partial\in\mathfrak g$, there is a unique vector
field $\partial^X $ such that
for any open subset $U$ of $X$ and any $f\in\mathcal O_X(U)$,
$\partial^X (f) (x)= \partial_g f(g\cdot x)|_{g=\mathbf 1}$.
The map $\partial\mapsto \partial^X $ is a section of
the restriction map.

\begin {lemm} \label{lemm.vectorfields}
For $G=\ga^n$ and for any $\partial\in\mathfrak g$ the restriction 
$\partial^X |_G$
is invariant under $G$.
\end{lemm}

(This is of course peculiar to the additive group.
For a general algebraic group~$G$, vector fields on~$G$ of
the form $\partial^X |_G$ need not be invariant.)

\begin{prop} \label{prop.logvanishing}
Let $x\in D$ and fix a local equation $s_D$ of $D$ in a neighborhood
$U$ of $x$.
Then,
for any $\partial\in\mathfrak g$, $\partial^X  \log s_D
= \frac{\partial^X (s_D)}{s_D} $
is a regular function in $U$.
\end{prop}
\begin{proof}
By purity, it is sufficient to prove this when $D$ is smooth
in a neighborhood of~$x$ (since $X$ is smooth,
a function which is regular in the complement to a codimension 2
subscheme is regular everywhere).
Then, we can choose ({\'e}tale) local coordinates
$x_1,\dots,x_n$ in $U$ (i.e.,~elements of $\mathcal O_X(U)$,
$\Omega^1_X|_U$ is free over $\mathcal O_U$ with basis $dx_1,\dots,dx_n$)
such that $s_D=x_1$.
Moreover, we can write uniquely
\[ \partial^X  = \sum_{j=1}^n f_j \frac{\partial}{\partial x_j}
\]
for some functions $f_j\in\mathcal O_X(U)$.
Let $D_x$ be the irreducible component of $D$ containing $x$.
Necessarily, $G$ stabilizes $D_x$ so that 
the function $g\mapsto x_1(g\cdot x)$ is identically
$0$ in a Zariski neighborhood of $\mathbf 1$ in $ G$ and a
fortiori, $\partial^X (x_1)$ vanishes on $U$.
By definition $\partial^X (x_1)=f_1$
hence $ f_1$ is a multiple of $x_1$: there exists a unique
$g_1\in\mathcal O_X(U)$ such that $f_1 = x_1 g_1$
and
\[ \partial^X  \log x_1 = g_1 \in\mathcal O_X(U). \]
The lemma is proved.
\end{proof}

\begin{exam}
For $G=\ga=\Spec F[x]$ or $G=\gm=\Spec F[x,x^{-1}]$, the Lie algebra 
$\mathfrak g$ has a canonical basis $\partial$, given
by the local parameter $x$ at the neutral element (respectively $0$ and $1$)
of $G$.
If we embed $G$ in $\P^1$ equivariantly, we get
$\partial^X =\partial /\partial x$ for $G=\ga$ and
$x \partial/\partial x$ for $G=\gm$. We see that it vanishes
at infinity (being $\{\infty\}$ or $\{0;\infty\}$, accordingly).
This is a general fact, as the following lemma shows.
\end{exam}

\begin{lemm}[Hassett/Tschinkel] \label{lemm.rho}
There exist integers $\rho_\alpha\geq 1$
such that 
\[ \omega_X^{-1}\simeq \mathcal
O_X(\sum \rho_\alpha D_\alpha). \]
If $G=\ga^n$, then for each $\alpha$,
$\rho_\alpha\geq 2$.
\end{lemm}
\begin{proof}
We only prove that $\rho_\alpha\geq 1$, because
the sharper bound which is valid in the case of $G=\ga^n$
won't be used below.
We refer to~\cite{hassett-t99}, Theorem 2.7 for its proof.

Let $\partial_1,\dots,\partial_n$ be a basis of $\mathfrak g$.
Then, $\delta:=\partial_1^X \wedge\cdots\wedge\partial_n^X $
is a global section of the line bundle $\det\mathcal T_X=\omega_X^{-1}$.
Moreover, $\delta$ does not vanish on $G$. Therefore,  we can write
$\div(\delta)=\sum\rho_\alpha D_\alpha$ for nonnegative
integers $\rho_\alpha$. Necessarily, $\omega_X^{-1}\simeq
\mathcal O_X(\sum\rho_\alpha D_\alpha)$, hence we have
to prove that these integers are positive.

Fix any $x$ in the smooth part of $X\setminus G$,
so that there exists a unique $\alpha\in\mathcal A$ such that  $x\in D_\alpha$.
Pick local coordinates $x_1,\dots,x_n$ in a neighborhood
$U$ of $x$ in such a way that in $U$, $D$ is defined by
the equation $x_1=0$.
Write $\partial_i^X =\sum f_{ij} \frac{\partial}{\partial x_j}$,
with $f_{ij}\in \mathcal O_X(U)$.
We have seen in Proposition~\ref{prop.logvanishing} that
for any $i$, $ \partial_i^X (x_1)=f_{i1}\in(x_1)$.
Hence, $\delta \in (x_1) \frac{\partial}{\partial x_1}\wedge\cdots
\wedge\frac{\partial}{\partial x_n}$.
\end{proof}


\section{Metrizations} \label{section.metrizations}

Let $F$ be a number field and $\mathfrak o_F$
its ring of integers. 
Denote by $F_v$ the completion of $F$ at a place $v$,
by  $\mathfrak o_v$ the ring of integers in $F_v$ 
if $v$ is nonarchimedean, by 
$\A_F$ the ring of adeles of $F$
and by $\Afin$ the (restricted) product of $F_v$ 
over the nonarchimedean places $v$.  
The valuation in $F_v$ is normalized in such a way 
that for any Haar measure $\mu_v$ on $F_v$ and any measurable
subset $I\subset F_v$, $\mu_v(a I)=\abs{a}_v \mu_v(I)$.
In particular, it is the usual absolute value for $F_v=\R$, its square
for $F_v=\C$, it satisfies $\abs{p}_p=1/p$ if $F_v=\Q_p$
and if $\mathcal N_v : F_v\ra \Q_p$ is the norm map, 
$\abs{x}_v=\abs{\mathcal N_v(x)}_p$.
Let $X$ be an equivariant compactification of $\ga^n$ as above and 
$\mathcal L$ a line bundle on $X$, endowed with its canonical
linearization.

\begin{defi}\label{defi.adelic}
A \emph{smooth adelic metric} on $\mathcal L$ is a family of $v$-adic norms
$\norm{\cdot}_v$ on $\mathcal L$ for all places $v$ of $F$ satisfying
the following properties:
\begin{enumerate}\def\theenumi{\alph{enumi}}\def\labelenumi{(\theenumi)}
\item if $v$ is archimedean, then $\norm{\cdot}_v$ is $\mathcal C^\infty$;
\item if $v$ is nonarchimedean, then $\norm{\cdot}_v$ is locally constant;
\par\noindent
\emph{(i.e., the norm of any local nonvanishing section is $\mathcal C^\infty$,
resp.\ locally constant) }
\item there exists an open dense subset $U\subset\Spec(\mathfrak o_F)$,
a flat projective $U$-scheme $\mathcal X_{/U}$ extending $X$
together with an action of $\ga^n{}_{/U}$ extending the action
of $\ga^n$ on $X$ and 
a linearized line bundle $\mathcal L$ on $\mathcal X_{/U}$
extending the linearized line bundle on $X$,
such that for any place $v$  lying over $U$,
the $v$-adic metric on $\mathcal L$ is given by the integral model.
\end{enumerate}
\end{defi}

\begin{lemm}
Let $v$ be a nonarchimedean valuation of $F$ and $\norm{\cdot}_v$
a locally constant $v$-adic norm on $\mathcal L$.
Then the stabilizer of $(\mathcal L,\norm{\cdot}_v)$,
i.e., the set of $g\in\ga^n(\mathfrak o_v)$ which 
act isometrically on $(\mathcal L,\norm{\cdot}_v)$,
is a compact open subgroup of $\ga^n(\mathfrak o_v)$.
\end{lemm}
\begin{proof}
First we assume that $\mathcal L$ is effective.
By Proposition~\ref{prop.invsection},
we have a nonzero invariant global section 
$\mathsec s$.

If $m$ and $p_2$ denote the action and the 
second projection $G\times X\ra X$, respectively,
endow the trivial line bundle
$m^*\mathcal L\otimes p_2^*\mathcal L^{-1}$ on $\ga^n\times X$
with the tensor-product metric. This is a locally constant
metric on the trivial line bundle.
The function
on $\ga^n(F_v)\times \ga^n(F_v)$ given by
\[
(g,x)\mapsto \norm{\mathsec s(g+x)}_v\norm{\mathsec s(x)}_v^{-1}.
\]
is the norm of the canonical basis $1$ and
therefore extends to a locally constant function
on $\ga^n(F_v)\times X(F_v)$.
Its restriction to the compact subset $\ga^n(\mathfrak o_v)\times X(F_v)$
is uniformly continuous. Since it is locally constant and
equal to $1$ on $\{1\}\times X(F_v)$, there
exists a neighborhood of $\{1\}\times X(F_v)$ on which it equals~$1$.
Such a neighborhood contains a neighborhood of the form $K_v \times X(F_v)$, 
where $K_v$ is a compact open subgroup of $\ga^n(\mathfrak o_v)$.
This proves the lemma in the effective case.

In the general case, we write
$\mathcal L =\mathcal L_1\otimes\mathcal L_2^{-1}$
for two effective line bundles (each having a nonzero global section).
We can endow $\mathcal L_2$ with any locally constant $v$-adic metric
and $\mathcal L_1$ with the unique (necessarily locally constant)
$v$-adic metric such
that the isomorphism $\mathcal L_1\simeq \mathcal L\otimes\mathcal L_2$
is an isometry.
By the previous case, 
there exist two compact open subgroups $K_{1,v}$ and $K_{2,v}$ 
contained in $\ga^n(\mathfrak o_v)$ which act
isometrically on $\mathcal L_1$ and $\mathcal L_2$, respectively.
Their intersection acts isometrically on $\mathcal L$.
\end{proof}

\begin{prop}\label{prop.stab}
If $\mathcal L$ is endowed with a smooth adelic metric then for all
but finitely many places $v$ of $F$
the stabilizer of $(\mathcal L,\norm{\cdot}_v)$
is equal to $\ga^n(\mathfrak o_v)$.
Therefore, their product over all 
finite places of $F$ is 
a compact open subgroup of $\ga^n(\Afin)$.
\end{prop}
\begin{proof}
It suffices to note
that if $v$ lies over the open subset $U\subset\Spec(\mathfrak o_F)$
given by the definition of an adelic metric then
the stabilizer of $(\mathcal L,\norm{\cdot}_v)$ equals
$\ga^n(\mathfrak o_v)$.
\end{proof}

{}From now on we choose
adelic metrics on the line bundles
$\mathcal L$ in such a way that
the tensor product of two line bundles is endowed with the
product of the metrics (this can be done by fixing smooth
adelic metrics on a $\Z$-basis of $\Pic(X)$ and extending 
by linearity).
We shall denote by $\mathbf K$ 
the compact open subgroup
of $\ga^n(\Afin)$ 
stabilizing all these metrized line bundles on $X$ .
We also fix
an open dense $U\subset \Spec(\mathfrak o_F)$,
a  flat and projective model $\mathcal X_{/U}$ 
over $U$ and models of the line bundles $\mathcal L$
such that the chosen $v$-adic metrics for all line bundles on $X$ 
are given by these integral models. 

\section{Heights} \label{section.heights}
 
Let $X$ be an algebraic variety 
over $F$ and $({\mathcal L}, \norm{ \cdot }_v)$
an adelically metrized line bundle on $X$. 
The associated  {\em height function} is defined as
\[
H_{\mathcal L} \,:\, X(F)\ra \R_{>0} 
\qquad
H_{\mathcal L}(x)=
\prod_v H_{{\mathcal L},v}(x):=\prod_v\norm{{\mathsec s}}_v(x)^{-1}, 
\] 
where ${\mathsec s}$ is any $F$-rational section of ${\mathcal L}$
not vanishing at $x$. The product formula ensures that
$H_{\mathcal L}$ does not depend on the choice of the 
$F$-rational section $\mathsec s$ (though the local heights
$H_{{\mathcal L},v}(x)$ do).

In Section~\ref{section.metrizations} we have defined 
simultaneous metrizations of line bundles on 
equivariant compactifications of $\ga^n$.  
This allows to define compatible systems of heights
\[ H\colon X(F) \times \Pic(X)_\C \ra \C . \]
Fix for any $\alpha\in\mathcal A$ some non zero $\ga^n$-invariant
section $\mathsec s_\alpha$ of $\mathscr O_X(D_\alpha)$.
We can then extend the height pairing $H$ to a pairing
\begin{equation}\label{eq.heightpairing}
H=\prod_v H_v \, :\, \ga^n(\A_F) 
\times \Pic(X)_{\mathbf C} \ra \C
\end{equation}
by  mapping
\[
({\mathbf x};{\mathbf s})=((\mathbf x_v);\sum s_\alpha D_\alpha))
   \mapsto 
   \prod_v \prod_\alpha \norm{\mathsec s_\alpha}_v(\mathbf x_v)^{-s_\alpha}.
\]

Recall that $\mathbf K$ denotes the compact open subgroup
of $\ga^n(\A_F)$ stabilizing all line bundles  on $X$ together
with their chosen metrization.
By construction, we have the following proposition.
\begin{prop}\label{prop.heights}
The height pairing $H$ defined in~\eqref{eq.heightpairing}
is $\mathbf K$-invariant in the first component
and (exponentially) linear
in the second component.
\end{prop}

\begin{prop}\label{prop.est}
Assume that $\mathcal L$ is in the interior of the effective
cone of $\Pic(X)$ (i.e., 
$\mathcal L\simeq \mathcal O_X(\sum d_\alpha D_\alpha)$
for some positive integers $d_\alpha>0$).
Then, for any real $B$, there are only finitely many $x\in\ga^n(F)$
such that $H(x;\mathcal L)\leq B$.
\end{prop}
In view of Lemma~\ref{lemm.rho}, this applies to the anticanonical
line bundle $K_X^{-1}$.
\begin{proof}
Let $\mathcal M$ be an ample line bundle on~$X$ and let $\nu$ be a sufficiently
large integer such that $\mathcal L^\nu\otimes\mathcal M^{-1}$ is effective.
It follows from the preceding section that
$\mathcal L^\nu\otimes \mathcal M^{-1}$ has a section $\mathsec s$
which does not vanish on $\ga^n$. This implies that
the function $x\mapsto -\log H(x;\mathcal L^\nu\otimes\mathcal M^{-1})$
is bounded from above on $\ga^n(F)$.
Therefore, there exists a constant $C>0$ such that for any $x\in\ga^n(F)$,
\[ H(x;\mathcal L) \geq C H(x;\mathcal M)^{1/\nu}.\]
We may now apply Northcott's theorem and obtain the
desired finiteness.
\end{proof}

\begin{rema}
The same argument shows that the rational map given
by the sections of a sufficiently high power of $\mathcal L$
is an embedding on $\ga^n$.
\end{rema}

The main tool in the study of asymptotics for the number
of points of bounded height is the
\emph{height zeta function} 
\[
Z (\mathbf s) = \sum_{x\in\ga^n(F)}
              H(x;{\mathbf s})^{-1} ,
\qquad \mathbf s=(s_\alpha)\in\Pic(X)_{\C}. \]

\begin{prop}\label{prop.abs-conv}
There exists an non-empty open subset $\Omega\subset\Pic(X)_\R$
such that 
$Z(\mathbf s)$ converges absolutely to a bounded holomorphic function in 
the tube domain $\Omega+i\Pic(X)_{\R}$ in
the complex vector space $\Pic(X)_{\C}$.
\end{prop}
\begin{proof}
Fix a basis $(\mathcal L_j)_j$ of $\Pic(X)$ consisting of 
(classes of) line
bundles lying in the interior of the effective cone of $\Pic(X)$.
Let $\Omega_t$ denote the open set of all linear combinations
$\sum t_j \mathcal L_j\in\Pic(X)_\R$ such that for some $j$,
$t_j> t$.
Fix some ample line bundle~$\mathcal M$. It is well known that
the height zeta function of $X$ relative to $\mathcal M$
converges for $\Re(s)$ big enough, say $\Re(s)> \sigma_0$. In the proof
of Proposition~\ref{prop.est} we may choose some $\nu$ which works
for any $\mathcal L_j$, so that for any $j$
\[ H(x;\mathcal L_j)^{-1} \ll H(X;\mathcal M)^{-1/\nu}. \]
Since $H(\cdot;\mathcal L_j)$ is bounded from below on $\ga^n(F)$
it follows that for any $\mathcal L\in \Omega_1$ 
\[ H(x;\mathcal L)^{-1} \ll H(X;\mathcal M)^{-1/\nu} \]
Therefore, the height zeta function converges absolutely and
uniformly on the tube domain $\Omega_{\nu\sigma_0}+i\Pic(X)_\R$.
\end{proof}
The following sections are devoted to the study of analytic
properties of $Z (\mathbf s)$.


\section{The Poisson formula} \label{section.poisson}

We recall basic facts concerning harmonic
analysis on the group $\G_a^n$ over the adeles $\A=\A_F$
(cf., for example, \cite{tate67b}).
For any prime number $p$,
we can view $\Q_p/\Z_p$ as the $p$-Sylow subgroup
of $\Q/\Z$ and we can define a local character $\psi_p$
of $\G_a(\Q_p)$
by setting
\[
\psi_p \colon  x_p \mapsto \exp(2 \pi i x_p).
\]
At the infinite place of $\Q$ we put
\[
\psi_{\infty}\colon  x_{\infty} \mapsto \exp(-2 \pi i x_{\infty}), 
\]
(here $x_{\infty}$ is viewed as an element in $\R/\Z$).
The product of local characters gives a character
$\psi$ of $\G_a(\A_{\Q})$ and, by composition with the trace,
a character of $\G_a(\A_F)$.

If $\mathbf a\in\ga^n(\A_F)$, 
recall that $f_{\mathbf a}=\langle \cdot,\mathbf a\rangle$
is the corresponding linear form on $\ga^n(\A_F)$ and let
$\psi_{\mathbf a}=\psi\circ f_{\mathbf a}$.
This defines a map $\ga^n(\A_F)\ra\ga^n(\A_F)^*$.  It
is well known that this map is a Pontryagin
duality. 
The subgroup $\ga^n(F)\subset \ga^n(\A_F)$ is discrete, cocompact and
we have an induced Pontryagin duality   
\[
\ga^n(\A_F)\ra (\ga^n(\A_F)/\ga^n(F))^*.
\] 
We fix selfdual Haar measures  $\mathrm dx_v$ on 
$\ga(F_v)$ for all $v$. We refer to \cite{tate67b}
for an explicit normalization of these measures. 
We will use the fact that for all but finitely many 
$v$ the volume of
$\ga(\mathfrak o_v)$ with respect to $\mathrm dx_v$ is equal to 1. 
Thus we have an induced selfdual Haar measure $\mathrm dx$ on $\ga(\A_F)$ and 
the product measure $\dx$ on $\ga^n(\A_F)$.
The Fourier transform (in the adelic component) 
of the height pairing on 
$\ga^n(\A_F)\times\Pic(X)_{\C}$ is defined by
\[ \hat H(\psi_{\mathbf a};\mathbf s)
        = \int_{\ga^n(\A_F)} H(\mathbf x;\mathbf s)^{-1} 
\psi_{\mathbf a}(\mathbf x)\, \dx .
\]
We will use the Poisson formula in the following form
(cf.~\cite{manin-p95}, p.~280).
\begin{theo}
\label{theo.poisson}
Let $\Phi$ be a continuous
function on $\G_a^n(\A_F)$ such that
both $\Phi$ and its Fourier transform $\hat{\Phi}$ are integrable
and such that the series
\[
\sum_{\mathbf x\in \G_a^n(F)} \Phi({\mathbf x}+\mathbf b )
\]
converges absolutely and uniformly when $\mathbf b$
belongs to $\G_a^n(\A_F)/\G_a^n(F)$.
Then, 
\[
\sum_{{\mathbf x}\in \G_a^n(F)} \Phi({\mathbf x})=
\sum_{ {\mathbf a}\in \G_a^n(F)}
\hat{\Phi}(\psi_{\mathbf a}).
\]
\end{theo}

The following lemma 
(a slight strengthening of Proposition~\ref{prop.abs-conv}) 
verifies the two hypotheses
of the Poisson formula~\ref{theo.poisson} concerning $H$. 

\begin{lemm} 
Let $X$ be a smooth projective
equivariant compactification of $\ga^n$ and $H$ the 
height pairing defined in  Section~\ref{section.heights}.
There exists a nonempty open 
subset $\Omega \subset \Pic(X)_{\R}$ such that
for any $\mathbf s\in \Omega +i\Pic(X)_{\R}$
the series
\[ \sum_{x\in\ga^n(F)} H(x+\mathbf b; \mathbf s)^{-1} \]
converges absolutely, uniformly in~$\mathbf b\in \ga^n(\A_F)/\ga^n(F)$
and locally uniformly in~$\mathbf s$.
Moreover, for such~$\mathbf s$, 
the function $H(\cdot;\mathbf s)^{-1}$ is integrable
over $\ga^n(\A_F)$.
\end{lemm}
\begin{proof}
Since the natural action of $\ga^n$ on $\Pic(X)$ is trivial,
for any $\mathbf b\in\ga^n(\A_F)$
the function $x\mapsto H(x+\mathbf b;\mathcal L)$ is a height function
for $\mathcal L$, induced by a ``twisted'' adelic metric. When
$\mathbf b$ belongs to some compact subset $K$ of $\ga^n(\A_F)$
there exists a constant $C_K(\mathcal L)$ such that for any $\mathbf
x\in\ga^n(\A_F)$
and any $\mathbf b\in K$,
\[   C_K(\mathcal L) ^{-1} H(\mathbf x;\mathcal L)
\leq H(\mathbf x+\mathbf b;\mathcal L)
\leq  C_K(\mathcal L)  H(\mathbf x;\mathcal L). \]
Indeed, the quotient $H(\mathbf x+\mathbf b;\mathcal L)/H(\mathbf x;\mathcal L)$
defines a bounded continuous function of
$(\mathbf x,\mathbf b)\in\ga^n(\A_F)\times K$.
(Only at a finite number of  places do they differ,
and at these places~$v$,
a comparison is provided by the compactness of~$X(F_v)$.)
As $\ga^n(\A_F)/\ga^n(F)$ is compact, one may take for $K$
any compact set containing a fundamental domain.
This implies the uniform convergence of the series
once $\mathbf s$ belongs to the  tube domain of Proposition~\ref{prop.abs-conv}.

The integrability follows at once : for such~$\mathbf s$,
\begin{align*}
 \int_{\ga^n(\A_F) } \abs{ H(\mathbf x;\mathbf s)}^{-1}\,\mathrm d\mathbf
x
& \leq \sum_{x\in\ga^n(F)} \int_{x+K} \abs{H(x+\mathbf b;\mathbf
s)}^{-1}\,\mathrm d\mathbf b \\
& \leq C_K \vol(K) \sum_{x\in\ga^n(F) } \abs{H(x;\mathbf s)}^{-1}.
\end{align*}
\end{proof}

The following proposition follows 
from the invariance of the height under the
action of the compact~$\mathbf K$.
\begin{prop}\label{prop.vanishing}
For all characters
$\psi_{\mathbf a}$ of $\ga^n(\A_F)$ which are nontrivial on the
compact subgroup $\mathbf K$ of $\ga^n(\Afin)$
and all $\mathbf s$ such that $H(\cdot ; \mathbf s)^{-1}$ is
integrable,
we have
\[ \hat{H}({\mathbf s};\psi_{\mathbf a})=0.\]
\end{prop}

Consequently, the 
set $\mathfrak d_X$ of all ${\mathbf a}\in\ga^n(F)$
such that $\psi_{\mathbf a}$ is 
trivial on ${\mathbf K}$ is is a sub-$\mathfrak o_F$-module
of $\ga^n(F)$,
commensurable with $\ga^n(\mathfrak o_F)$.

Consequently, provided the Poisson formula applies
and we have a  formal identity for the
height zeta function:
\begin{equation} \label{eq.formal}
Z({\mathbf s})= 
   \sum_{\mathbf a\in\mathfrak d_X} 
    \hat{H}({\mathbf s};\psi_{\mathbf a}).
\end{equation}

The last four Sections of the paper are concerned with
analytic arguments leading to
the evaluation of $\hat H$ and to the fact that Equation~\eqref{eq.formal}
holds for all $\mathbf s$ in some tube domain.
It will be necessary 
to assume that the boundary divisor~$D$ has strict normal crossings,
which means that over the separable closure~$\bar F$ of~$F$,
$D$ is a union of smooth irreducible components meeting transversally.

We sum up these Sections in the following Proposition.
Recall that for $\mathbf a\in\ga^n(F)$, $\mathcal A_0(\mathbf a)$
is the set of $\alpha\in\mathcal A$  such that $d_\alpha(f_{\mathbf a})=0$.
In particular, $\mathcal A_0(0)=\mathcal A$.
Also, for any $t\in\R$, we denote by $\Omega_t$ the open tube
domain in $\Pic(X)_\C$ defined by the inequalities $s_\alpha>\rho_\alpha+t$.
Finally, let $\norm{\cdot}_\infty$ denote any norm
on the real vector space
$\ga^n(F\otimes_\Q\R)\simeq \prod_{v|\infty} F_\infty^n$.
\begin{prop}\label{prop.main-calc}
Assume that $D$ has strict normal crossings.

\textup{a)}
For any $\mathbf s\in\Omega_0$,
$H(\cdot;\mathbf s)$ is integrable on $\ga^n(\A_F)$.

\textup{b)}
For any $\eps >0$ and $\mathbf a\in \mathfrak d_X $
there exists a holomorphic bounded function $\phi(\mathbf a;\cdot)$
on $\Omega_{-1/2+\eps}$ such
that for any $\mathbf s\in \Omega_{0}$ 
\[ \hat H(\psi_{\mathbf a};\mathbf s)=
\phi(\mathbf a;\mathbf s)\prod_{\alpha\in \mathcal A_0(\mathbf a)}
     (s_\alpha-\rho_\alpha)^{-1}. \]

\textup{c)}
For any $N>0$ there exist  
constants $N'>0$ and $C(\eps, N)$ such that
for any $\mathbf s\in\Omega_{-1/2+\eps}$
and any $\mathbf a\in\mathfrak d_X$
one has the estimate
\[ \abs{\phi(\mathbf a;\mathbf s)} \leq C(\eps,N) 
     (1+\norm{\Im(\mathbf s )})^{N'}
     (1+\norm{\mathbf a}_\infty)^{-N}. \]
\end{prop}

\section{Meromorphic continuation of the height zeta function}
\label{sec.singularities}

For the proof of our main theorem,
we need to extend the estimates of Proposition~\ref{prop.main-calc}
when $D$ is not assumed to have strict normal crossings.

\begin{lemm}\label{lemma.singularities}
Proposition~\ref{prop.main-calc} remains true without the
hypothesis that $D$ has strict normal crossings.
\end{lemm}

\begin{proof}
Equivariant resolution of singularities (in char.~$0$,
see~\cite{encinas-v98} or~\cite{bierstone-m97})
shows that there exists a 
composition of equivariant blow-up with smooth centers $\pi\colon\tilde
X\ra X$ such that 
\begin{itemize}
\item $\tilde X$ is a smooth projective equivariant compactification
of $\ga^n$;
\item $\pi$ is equivariant (hence induces an isomorphism on the
$\ga^n$);
\item the boundary divisor $\tilde D=\tilde X\setminus \ga^n$
has strict normal crossings.
\end{itemize}
The map $\pi$ induces a map from the
set $\tilde {\mathcal A}$ of irreducible components of $\tilde D$
to~$\mathcal A$. Moreover, $\pi$
has a canonical section obtained by associating
to a component of $D$ its strict transform in $\tilde X$.
This allows to view $\mathcal A$ as a subset of $\tilde {\mathcal A}$,
the complementary
subset $\mathcal E$ consisting of exceptional divisors.
Via the identifications of $\Pic(X)$ with $\Z^{\mathcal A}$
and $\Pic\tilde X$ with $\Z^{\tilde{\mathcal A}}$,
the map $\pi^*$ on line bundles
induces a linear map $i\colon \Z^{\mathcal A} \ra  \Z^{\tilde{\mathcal A}}$
which is a section of the first projection 
$\Z^{\tilde{\mathcal A}}=\Z^{\mathcal A}\oplus\Z^{\mathcal E}
\ra \Z^{\mathcal A}$.  (We also extend $i$ by linearity to the vector
spaces of complexified line bundles.)

Concerning metrizations, it is possible to endow line bundles
on $\tilde X$ with adelic metrics which extend the metrics
coming from $X$ via $\pi^*$. Remark that with such a choice,
the height function on $\ga^n(\A_F)$ and its Fourier transform
on $\tilde X$ extend those on $X$.
For instance,
for any $\mathbf s\in\Pic(X)_\C$ and any $\mathbf x\in\ga^n(\A_F)$,
one has
\[ H(\mathbf x;\mathbf s) = \tilde H(\mathbf x;i(\mathbf s)). \]
As $i$ maps $\N^{*\mathcal A}$ into $\N^{*\tilde{\mathcal A}}$,
it follows in particular
that $H(\cdot;\mathbf s)$ is integrable for  $\mathbf s\in\Omega_0$
and that for any $\mathbf a\in\ga^n(\A_F)$ and any $\mathbf s\in\Omega_0$,
one has
\[
\hat H(\psi_{\mathbf a};\mathbf s )
        = \Hat{\tilde H}(\psi_{\mathbf a}; i( \mathbf s)). \]

Unless $\pi$ is an isomorphism, $\pi$ is not smooth
and $\omega_{\tilde X}\otimes\pi^*\omega_X^{-1}$ has a nonzero
coefficient at each exceptional divisor.
For any effective $\lambda$,
this implies that $i(\lambda)-a_{\lambda}\tilde\rho$
lies on the boundary of $\Lambda_\eff(\tilde X)$.
(Recall that the classes of $\omega_X^{-1}$ and $\omega_{\tilde X}^{-1}$
in $\Z^{\mathcal A}$ and $\Z^{\tilde{\mathcal A}}$ are denoted
$\rho$ and $\tilde\rho$.)
It follows that
$ a_{i(\lambda)}=a_{\lambda}$ and  that
\[ {\mathcal B}(i(\lambda))
  = \{ \tilde\alpha\in\tilde{\mathcal A}\,;\,
   i(\lambda) =  a_{i(\lambda)}\tilde\rho \}
\]
can be identified with ${\mathcal B}(\lambda)$.

Consequently, any ``exceptional'' factor
$(s_\alpha-\tilde\rho_\alpha)^{-1}$ for $\alpha\in\mathcal E$
extends to a holomorphic bounded function on $i(\Omega_{-1/2})$
This shows that Proposition~\ref{prop.main-calc} still stands true
without the hypothesis that $D$ has strict normal crossings.
\end{proof}

One has also the following result which shall allow
for the determination of the constant in the asymptotics.
\begin{prop}\label{prop.tamagawa}
Denote by $\tau(K_X)$ the Tamagawa number of $X$
defined by the adelic metrization on $K_X$. One then
has the equality
\[ \phi(0;\rho)=\tau(K_X)  . \]
In particular, it is a strictly positive real number.
\end{prop}
\begin{proof}
Recall first Peyre's definition~\cite{peyre95}
of the Tamagawa measure
in this context, which extends Weil's definition in the context
of linear algebraic groups~\cite{weil82}.
For any place $v$ of $F$ and 
any local non-vanishing differential form $\omega$ on $X(F_v)$,
the self-dual measure on $\ga^n(F_v)$ induces
a local measure $\abs{\omega}/\norm{\omega}$ on $X(F_v)$. These
local measures glue to define a global measure $\mathrm d\mu_v$
on $X(F_v)$.

Recall that Artin's $L$-function of $\Pic(X_{\bar F})$
is equal to $\prod_\alpha \zeta_{F_\alpha}$.
Using Deligne's theorem on Weil's conjectures
(but it will also follow from our calculations),
Peyre has shown
that the Euler product
\[ \prod_v \left( \mu_v(X(F_v)) \prod_\alpha
\zeta_{F_\alpha,v}(1)^{-1}\right) \]
converges absolutely and he defined the Tamagawa number of $X$
as
\[ \tau(K_X)= \prod_\alpha \res_{s=1}\zeta_{F_\alpha}(s)
 \times \prod_v \left( \mu_v(X(F_v)) \prod_\alpha
\zeta_{F_\alpha,v}(1)^{-1}\right).
\]
(Note however that Peyre
apparently doesn't use the selfdual measure
for his definition in \emph{loc. cit.},
but inserts the appropriate correcting factor $\Delta_F^{-\dim X/2}$.)

On the open subset $\ga^n(F_v)$,
we can define the local measure
with the differential form $\omega=\mathrm dx_1
\wedge\dots\wedge \mathrm dx_n$ induced by the canonical
coordinates of $\ga^n$. Let $c\in F^\times$ be such that
$\omega= c\prod_\alpha {\mathsec s_\alpha}^{-\rho_\alpha}$.
On $\ga^n(F_v)$, we thus have
\[ d\mu_v = \abs{c}_v^{-1} \prod_\alpha \norm{\mathsec
s_\alpha}(x)^{-\rho_\alpha}\, \mathrm dx
     = \abs{c}_v^{-1} H_v(x;\rho)\;\mathrm dx \]
so that, the boundary subset being of measure~$0$,
\[ \mu_v(X(F_v))=\mu_v(\ga^n(F_v))=
\abs{c}_v^{-1}\int_{\ga^n(F_v)}H_v(x;\rho)\,\mathrm dx
      = \abs{c}_v^{-1}\hat H_v(\psi_0;\rho).\]

Moreover, it follows from these results and the absolute convergence
of the Euler products of Dedekind zeta functions for $\Re(s)>1$
that for any $s>1$, $H(\cdot;s\rho)$ is integrable on  $X(\A_F)$
and that one has
\begin{align*}
\phi(0;\rho) &= \lim_{s\ra 1^+} \phi(0;s\rho) \\
&= \lim_{s\ra 1^+} \hat H(\psi_0;s\rho) \prod_{\alpha\in\mathcal A}
        (s\rho_\alpha-\rho_\alpha)) \\
&= \prod_\alpha
\big(\lim_{s\ra 1^+}
\rho_\alpha(s-1)\zeta_{F_\alpha}(1+\rho_\alpha(s-1))
\big)
\times \\
& \hskip 6pc
\lim_{s\ra 1^+} \prod_\alpha \zeta_{F_\alpha}(1+\rho_\alpha(s-1))^{-1}
\prod_v \hat H_v(\psi_0;s\rho) \\
&=
 \prod_\alpha\res_{s=1}\zeta_{F_\alpha}(s)
\times
\lim_{s\ra 1^+} \prod_v \hat H_v(\psi_0;s\rho)\prod_\alpha\zeta_{F_\alpha,v}(1+\rho_\alpha(s-1))^{-1} 
\end{align*}
The estimates of Proposition~\ref{prop.bbb} below allow
us to interchange the limit $\lim_{s\ra 1}$ and the infinite product
$\prod_v$, so that
\begin{align*}
 \phi(0;\rho) & =
\prod_\alpha\res_{s=1} \zeta_{F_\alpha}(s)
\times
\prod_v \left(\hat H_v(\psi_0;\rho) \prod_\alpha
\zeta_{F_\alpha,v}(1)^{-1} \right) \\
&= 
  \tau(K_X) \prod_v \abs{c}_v = \tau(K_X),
\end{align*}
because of the product formula.
\end{proof}

It follows from these results that the Poisson formula~\eqref{eq.formal}
applies for any $\mathbf s\in\Omega_0$ and that
\[ 
 Z(\mathbf s) = \sum_{\mathbf a \in \mathfrak d_X} 
\hat H(\psi_{\mathbf a};\mathbf s) 
= \sum_{\mathbf a\in\mathfrak d_X} \phi(\mathbf a;\mathbf s)
  \prod_{\alpha\in \mathcal A_0(\mathbf a)}
        (s_\alpha-\rho_\alpha)^{-1}.
\]
Hence,
\[ 
Z(\mathbf s)\prod_{\alpha\in \mathcal A} (s_\alpha-\rho_\alpha)
 =\sum_{\mathbf a\in\mathfrak d_X }
      \phi(\mathbf a;\mathbf s)
       \prod_{\substack{ \alpha\in \mathcal A \\
                           \alpha\not\in \mathcal A_0(\mathbf a) }}
(s_\alpha-\rho_\alpha).
\] 
This last series is a sum of holomorphic functions on $\Omega_{-1/2}$
and taking $N>n[F:\Q]$, the estimate in Proposition~\ref{prop.main-calc}
implies that the series converges locally uniformly in that domain.
Therefore, we have shown that
the function
\[ \mathbf s\mapsto Z(\mathbf s)
\prod_{\alpha\in \mathcal A}(s_\alpha-\rho_\alpha)\]
has a holomorphic continuation to $\Omega_{-1/2}$
with polynomial growth in vertical strips.
The restriction of $Z(\mathbf s)$ to a line 
$\C(\lambda_\alpha)$ (with $\lambda_\alpha>0$ for all $\alpha$) 
gives 
\[
Z(s\lambda)= h_{\lambda}(s) \times 
\prod_{\alpha\in \mathcal A} (s\lambda_\alpha -\rho_\alpha) ^{-1}
\]
where $h_{\lambda}$ is a holomorphic function 
for $s\in\C$ such that
\[ \Re(s)>\max((\rho_\alpha-\frac12)/\lambda_\alpha)
     = a_{\lambda} - \frac{1}{2\min \lambda_\alpha},
\]
providing 
a meromorphic continuation of $Z(s\lambda)$ to this domain. 
The rightmost pole is at
\[  a_\lambda = \max_\alpha (\rho_\alpha/\lambda_\alpha)   \]
and its order is less or equal than $b_\lambda$.

\begin{theo}
Let $X$ be a smooth projective equivariant compactification
of $\ga^n$ over  a number field~$F$.
For any strictly effective $\lambda\in\Pic(X)$, the height
zeta function $Z(s\lambda)$ converges absolutely for $\Re(s)>a_\lambda$,
has a meromorphic extension to the left of $a_\lambda$,
with a single pole at $a_\lambda$ of order $\leq b_\lambda$.
\end{theo}

When $\lambda=\rho$, $a_\rho=1$, $b_\rho=\rk\Pic(X)=r$ and
Proposition~\ref{prop.tamagawa} implies the equality
\[ \lim_{s\ra 1} Z(s\rho)(s-1)^r
 = \tau(K_X)\prod_{\alpha\in \mathcal A}
\rho_\alpha^{-1} ,
\]
hence the order of the pole is exactly $r$. Using a standard Tauberian
theorem, we deduce an asymptotic expansion for the number
of rational points in $\ga^n$ of bounded anticanonical height.
\begin{theo}
\label{theo.rho}
Let $X$ be a smooth projective equivariant compactification of $\ga^n$
over a number field~$F$.

There exists a real number $\delta >0$ and a polynomial $P\in\R[X]$
of degree $r-1$
such that 
the number of $F$-rational points on $\ga^n\subset X$ 
of anticanonical height $\leq B$ satisfies
\[ N(K_X^{-1},B) = B P(\log B) + O(B^{1-\delta}). \]
The  leading coefficient of $P$ is given by
 \[ \frac{1}{(r-1)!}
    \tau(K_X)
    \prod_{\alpha\in \mathcal A} \rho_\alpha^{-1}. \]
\end{theo}

\section{Asymptotics}\label{section.general}

Let $\lambda\in \Lambda_\eff(X)$ be an effective class
and $\mathcal L_{\lambda}$ the corresponding line
bundle on $X$ equipped with a smooth adelic metric as in
Section~\ref{section.metrizations}.

The function $s\mapsto Z(s\lambda)$ is holomorphic
for $\Re(s)>a_{\lambda}=\max(\rho_\alpha/\lambda_\alpha)$ and
admits a meromorphic continuation to the left of $a_{\lambda}$.
Recall from  the previous section that it has a pole of order
$\leq b_{\lambda}$ at $a_{\lambda}$, where $b_{\lambda}$ is the 
cardinality of
\[ {\mathcal B}_{\lambda} 
= \{ \alpha\in\mathcal A\, ;\, \rho_\alpha=a_\lambda
\lambda_\alpha\}. \]
(Geometrically, the integer $b_{\lambda}$ is the codimension of
the face of the effective cone $\Lambda_\eff(X)$ containing
the class $a_{\lambda}\lambda-\rho$.)

In this section we prove that the order of the pole
of the height zeta function $Z(s\lambda)$ at $s=a_\lambda$
is exactly $b_\lambda$ and derive our main theorem.

Denote by $Z_\lambda(s)$ the sub-series
\[ Z_\lambda(s)= \sum_{\substack{\mathbf a \,\,\text{such that} \\
                    \mathcal A_0(\mathbf a)\supset {\mathcal B}_{\lambda}}}
             \hat H(\mathbf a;s\lambda). \]
It follows from the calculations above that
\begin{equation}
\lim_{s\ra a_{\lambda}} Z(s \lambda)(s-a_{\lambda})^{b_{\lambda}}
  = \lim_{s\ra a_{\lambda}} Z_\lambda(s) (s-a_{\lambda})^{b_{\lambda}} 
\end{equation}
since all other
terms converge to zero when $s\ra a_{\lambda}$,
with uniform convergence of the series.

The set~$V_\lambda$ of $\mathbf a\in \ga^n(F)$
such that $\mathcal A_0(\mathbf a)$ contains
${\mathcal B}_{\lambda}$ is a sub-vector space. 
Let $\mathbf G_\lambda\subset\ga^n$ be the $F$-sub-vector space defined
by the corresponding linear forms $\langle \mathbf a,\cdot\rangle$.
Then, the autoduality on $\ga^n(\A_F)$ identifies $V_\lambda$
with the dual of $\mathbf G_\lambda(\A_F)\ga^n(F)$.
We apply the Poisson summation formula and obtain
\[ Z_\lambda(s) = \int_{G_\lambda(\A_F)\ga^n(F)} 
       H(\mathbf x;s\lambda)^{-1}\, \dx
    = \sum_{\mathbf x\in(\ga^n/\mathbf G_\lambda)(F)} 
        \int_{\mathbf G_\lambda(\A_F)} 
H(\mathbf x+\mathbf y;s\lambda)^{-1}\,
\mathrm d\mathbf y.  \]
(The justification of Poisson summation formula is as
in Section~\ref{section.poisson}.)

\begin{prop}\label{prop.tamlambda}
For each $x\in (\ga^n/\mathbf G_\lambda)(F)$ there exists
a strictly positive real $\tau_\lambda(x)>0$ such
that
\[
\lim_{s\ra a_{\lambda}^+} (s-a_{\lambda})^{b_{\lambda}}
  \int_{\mathbf G_\lambda(\A_F)} H(x+\mathbf y;s\lambda)^{-1}\,\mathrm
d\mathbf y 
   =\tau_\lambda(x). \]
\end{prop}

\begin{proof}
It is sufficient to prove the proposition when $x=0$ as
the integrals for different values of $x$ are comparable.
The Zariski closure
$Y$ of $\mathbf G_\lambda$ in $X$ need not be smooth.
Nevertheless, we may introduce an equivariant proper modification
$\pi\colon\tilde X\ra X$ such that the Zariski closure $\tilde
Y\subset\tilde X$ of $\mathbf G_\lambda$ is a smooth equivariant
compactification whose boundary is a divisor with strict normal
crossings obtained by intersecting the components of the boundary
of $\tilde X$ with $\tilde Y$.
The arguments of Lemma~\ref{lemma.singularities}
show that we do not lose any generality by assuming
this on~$X$ itself.
\let\qedsymbol\relax
\end{proof}

\begin{lemm}
There exist integers $\rho'_\alpha\leq\rho_\alpha$ for 
$\alpha\in\mathcal A$
such that 
\[ \omega_{Y}^{-1} =\sum\rho'_\alpha (D_\alpha\cap Y). \]
Moreover, $\rho'_\alpha=\rho_\alpha$ if and only if
$\alpha\in \mathcal B(\lambda)$.
\end{lemm}
\begin{proof}
We prove this by induction on the codimension of $Y$ in $X$.
If $\mathbf G_\lambda=\div(f_{\mathbf a})$, the adjunction formula shows
that
\[ \omega_Y^{-1} = \omega_X^{-1}(-Y)|_Y = \sum
(\rho_\alpha-d_\alpha(\mathbf a))(D_\alpha\cap Y). \]
For $\mathbf a\in V_\lambda$ the nonnegative integer $d_\alpha(\mathbf a)$
is equal to~$0$ if and only if $\alpha\in {\mathcal B}_{\lambda}$.
The lemma is proved.
\end{proof}

\begin{proof}[End of the proof of Proposition~\ref{prop.tamlambda}]
Let $\mathrm d\tau'$ be the (nonrenormalized) Tamagawa measure
on $Y(\A_F)$: by definition, on $\mathbf G_\lambda(\A_F)$,
$\mathrm d\tau=H(\mathbf y;\rho')^{-1}\mathrm d\mathbf y$.
We have to estimate
\[ \lim_{s\ra a_{\lambda}} (s-a_{\lambda})^{b_{\lambda}} 
  \int_{Y(\A_F)} H(\mathbf y; s\lambda-\rho')^{-1}\,\mathrm d\tau'.\]
The integral in the limit is of the type studied in Prop~\ref{prop.main-calc},
but on the subgroup~$Y$.
Denote therefore by primes $'$ restrictions of objects from~$X$
to~$Y$.
Thanks to the previous lemma, one has $a_{\lambda'}=a_{\lambda}$,
$\mathcal B(\lambda')=\mathcal B(\lambda)$ and
and $b_{\lambda'}=b_\lambda$.
Therefore, with notations from Proposition~\ref{prop.main-calc},
\begin{align*}
 \lim_{s\ra a_{\lambda}} (s-a_{\lambda})^{b_{\lambda}} 
  \int_{Y(\A_F)} H(\mathbf y; s\lambda-\rho')^{-1}\,\mathrm d\tau' 
 \hskip -4cm \\
& =
 \lim_{s\ra a_{\lambda'}} (s-a_{\lambda'})^{b_{\lambda'}}
\prod_{\alpha\in \mathcal B(\lambda)}
\zeta_{F_\alpha}(1+s\lambda'_\alpha-\rho'_\alpha)
\phi'(s\lambda'-\rho',0)  \\
& = \phi'(0,0) >0 .
\end{align*}
The proposition is proved.
\end{proof}

We now can conclude as in the case of the anticanonical line bundle:
\begin{theo} \label{theo.lambda}
Let $\lambda\in \Pic(X)$ be the class of a line bundle 
contained in the interior of the effective cone $\Lambda_{\eff}(X)$,
and equip  the corresponding
line bundle  $\mathcal L_{\lambda}$
with a smooth adelic metric as in 
Section~\ref{section.metrizations}. 
There exist a polynomial $P_\lambda\in\R[X]$ of degree~$b_{\lambda}-1$
and a real number $\eps>0$
such that 
the number $N(\mathcal L_{\lambda}, B)$ 
of $F$-rational points on $\ga^n$ of $\mathcal L_{\lambda}$-height
$\leq B$ satisfies
\[
N(\mathcal L_{\lambda},B)=B^{a_{\lambda}} 
P_\lambda(\log B) + O(B^{a_{\lambda}-\eps}).
\]
\end{theo}

The leading term of $P_\lambda$ is equal to
\[ \frac{1}{(b_{\lambda}-1)!}
    \prod_{\alpha\in {\mathcal B}_{\lambda}} \lambda_\alpha^{-1}
     \left(
       \sum_{x\in(\ga^n/\mathbf G_\lambda)(F)}
         \tau_\lambda(x)
     \right)
.\]
It is compatible with the description by Batyrev and Tschinkel
in~\cite{batyrev-t98}.


\section{General estimates} \label{section.estimates}

In the remaining sections we have to prove Prop.~\ref{prop.main-calc}.
First we define a set of bad places $S$:
as in Section~\ref{section.metrizations}, fix 
a good model $\mathcal X_{/U}$ (flat projective $U$-scheme) 
over a dense open subset 
$U\subset \Spec(\mathfrak o_F)$.
In particular, we assume that $\mathcal X_{/U}$ extends $X$
as an equivariant compactificatation of $\ga^n$
whose boundary consists of the schematic closures of the $D_\alpha$.
Moreover, we can assume that for any $\alpha$,
the chosen section $\mathsec s_\alpha$
of $\mathscr O(D_\alpha)$ extends to $\mathcal X_{/U}$ and
that its norm coincides with the one that can be defined with
this model at all finite places dominating~$U$.
In particular, at such places~$v$,
the local height functions restricted to $\ga^n(F_v)$ 
are invariant under $\ga^n(\mathfrak o_v)$. 
We can also assume that for any finite place~$v$
dominating~$U$, the (nonrenormalized) Tamagawa measure
$\mathrm d\mu_v$ is given on $\ga^n(F_v)$
by $H_v(x;\rho)\,\mathrm d\mathbf x$.

Let $S$ be the set of places~$v$ of~$F$ such that either
\begin{itemize}
\item $v$ is archimedean, or 
\item $v$ doesn't dominate~$U$, or 
\item the residual characteristic of $v$ is $2$ or~$3$, or
\item the volume of $\ga^n(\mathfrak o_v)$ with respect to 
$\dx_v$ is not equal to 1 (equivalently, $F$ is ramified at~$v$), or
\item over the local ring $\mathfrak o_v$, the union
$\bigcup_\alpha D_\alpha$ is not a transverse union of
smooth relative divisors over~$U$.
\end{itemize}
For $\mathbf a\in\mathfrak d_X$,
let $f_{\mathbf a}$ be the corresponding linear form on $\ga^n$, considered
as an element of $F(X)$ or of $F(\mathcal X_{/U})$
and $\div(f_{\mathbf a})$ its
divisor.
We denote by $S({\mathbf a})$ the set
of all valuations~$v$ such that either
\begin{itemize}
\item $v$ is contained in $S$, or 
\item $\div(f_{\mathbf a})$ is not flat over~$U$,
\emph{i.e.}
if $\mathfrak m_v$ denotes the maximal ideal of~$v$ in $\mathfrak o_F$,
$\mathbf a$ belongs to $\mathfrak m_v\mathfrak d_X$.
\end{itemize}

Finally, recall that for any $t\in\R$, $\Omega_t$ is the tube
domain in $\Pic(X)_\C$ consisting of classes $(s_\alpha)$ such that
for any $\alpha$,
$\Re(s_\alpha)>\rho_\alpha$.

In Section~\ref{section.estimates}  
we prove estimates for finite products
$\prod_{v\in S({\mathbf a})}\hat{H}_v(\psi_{\mathbf a};{\mathbf s})$. 
In Sections~\ref{section.trivchar}
and~\ref{section.otherchars} 
we compute explicitly the local Fourier transforms 
$\hat{H}_v( \psi_{\mathbf a} ; {\mathbf s} )$
for all $v\notin S({\mathbf a})$. 

\emph{In the following four sections we temporarily assume that
the irreducible components of the boundary of $X$ are
geometrically irreducible.}
In Section~\ref{section.nonsplit},
we shall explain how our results extend to the general case,
(still assuming that $D$ has strict normal crossings).

\bigskip

\begin{prop}\label{prop.estim}
Let $\Sigma$ be any finite set of places of~$F$
containing the archimedean places.
For every $\epsilon >0$
and any $N>0$ there exists a constant $C(\Sigma,\epsilon,N)$ such that for  
all $\mathbf a\in \mathfrak d_X$ and 
all $\mathbf s =(s_{\alpha})\in \Omega_{-1+\eps}$
one has the estimate
\[
\abs{\prod_{v\in \Sigma} \hat{H}_v( \psi_{\mathbf a} ; \mathbf s )} 
\le C(\Sigma, \epsilon,N) 
\frac{(1+\norm{\mathbf s})^{N[F:\Q]}}{\prod_{v|\infty}(1+\norm{\mathbf a}_v)^N},
\]
where $\norm{\mathbf s}= \norm{\Re({\mathbf s})}+\norm{\Im({\mathbf s})}$.
\end{prop}

We subdivide the proof of this proposition into a sequence of 
lemmas.

\begin{lemm}\label{prop.integrability}
The function $H_v(\cdot ; \mathbf s)^{-1}$
is integrable on $\ga^n(F_v)$ if and only if 
$\mathbf s\in\Omega_{-1}$ (\emph{i.e.},
$\Re(s_\alpha)>\rho_\alpha-1$ for all  $\alpha\in\mathcal A$). 
Moreover, for all $\epsilon >0$ 
and all nonarchimedean places $v$,
there exists a constant $C_v(\epsilon)$ 
such that for all  
$\mathbf s\in \Omega_{-1+\epsilon}$
and all ${\mathbf a}\in \mathfrak d_X$ one has the estimate
\[
\abs{\hat{H}_v( \psi_{\mathbf a} ; \mathbf s )}\le C_v(\epsilon).
\] 
\end{lemm}
\begin{proof}
Without loss of generality, we can assume
that $\psi_{\mathbf a}(\mathbf x)\equiv 1$.
Now, using an analytic partition of unity on $X(F_v)$ and the
assumption that
the boundary $D=\sum_{\alpha\in \mathcal A}D_\alpha$ is
a divisor with strict normal crossings,
we see that it suffices to compute the integral over a
relatively compact neighborhood
of the origin in $F_v^n$ (denoted by $\mathcal B$) on which we have
a set coordinates  $x_1,\dots,x_n$ so that in $\mathcal B$
the divisor $D$ is defined by the
equations $x_1\cdots x_a=0$ for some $a\in\{1,\dots,n\}$.
In $\mathcal B$ there exist continuous
bounded functions $h_{\alpha}$ (for $\alpha\in \mathcal A$) such that
\[ \int_{\mathcal B} H_v(\mathbf x;\mathbf s)^{-1}\dx_v
   = \int_{\mathcal B} \prod_{i=1}^a 
\abs{x_i}_{v}^{s_{\alpha(i)}-\rho_{\alpha(i)}}
        \exp(\sum s_\alpha h_\alpha(x)) \, \mathrm d\mu_v .
\]
The integral over ${\mathcal B}$ is comparable to an integral
of the same type with functions $h_\alpha$ replaced by $0$.
The lemma is now a consequence of the following
well-known lemma.
\end{proof}

\begin{lemm}
Let $K$ be a local field.  
The function $x\mapsto \abs{x}^s$ is integrable on the unit
ball in $K$ if and only if $\Re(s)>-1$.
\end{lemm}
\begin{proof}
We may assume $s\in\R$.
Choose $c\in\left]0;1\right[$ such that the annulus $c<\abs{x}\leq 1$
has positive measure in $K$ and let
$ I_0=\int_{c<\abs{x}\leq 1} \abs{x}^s\, \mathrm dx$.
Then, we have
\[
 I_n  =\int_{c^{n+1}<\abs{x}\leq c^n} \abs{x}^s\mathrm dx
     = c^{n(s+1)} I_0. 
\]
It follows that the integral over $K$ converges if and only if the
geometric series $\sum c^{n(s+1)}$ converges, that is 
if $s+1>0$.
\end{proof}

\begin{prop}
Let $v$ be an archimedean place of $F$ and let $\eps_v=[F_v:\R]$.
For any $N>0$ and any $\eps>0$, there exists a constant
$C_v({\eps, N})$ such that for any 
$\mathbf s\in \Omega_{-1+\eps}$
and all $\mathbf a\in\mathfrak d_X$ (and $\neq 0$), we have
\[ \hat H_v( \psi_{\mathbf a} ; \mathbf s ) 
\leq C_v({\eps,N}) 
\frac{(1+\norm{\mathbf s})^{\eps_v N}}{(1+\norm{\mathbf a}_v)^N}. \]
\end{prop}
\begin{proof}
Let us assume for the moment that 
$F_v=\R$. We shall write $ \mathbf a_v = (a_1,...,a_n)$,
$\mathbf x_v=(x_1,...,x_n)$ and $\mathrm d\mathbf x_v = 
\mathrm d\mathbf x$.  
Then, in the domain $\Omega_{-1}$  
one has
\[ \hat H_v( \psi_{\mathbf a} ; \mathbf s )
  = \int_{\R^n} H_v(\mathbf x ; \mathbf s)^{-1} 
\exp(-2i\pi\langle\mathbf a,\mathbf x\rangle)
\, \dx . \]
Now integrate by parts: for  any
$j\in\{1,\dots,n\}$ we have
\[ 2i\pi a_j \hat H_v( \psi_{\mathbf a} ; \mathbf s )
   =  \int_{\R^n} \frac{\partial}{\partial x_j} 
       H_v(\mathbf x ;\mathbf s)^{-1}
        \exp(-2i\pi\langle\mathbf a,\mathbf x\rangle)\,\dx \]
and by induction
\[ (2i\pi a_j)^N \hat H_v( \psi_{\mathbf a} ; \mathbf s )
   = \int_{\R^n} 
     \left( \frac{\partial^N}{\partial x_j^N} 
      H_v(\cdot ;\mathbf s)^{-1} \right)(\mathbf x)
\exp(-2i\pi\langle\mathbf a,\mathbf x\rangle)\,\dx. \]
For any $\alpha$, let us define
$h_\alpha=\log\norm{\mathsec s_{D_\alpha}}$;
it is a $\mathcal C^\infty$ function on $\ga^n(\R)$.

Let $x$ be a point of $X(\R)$ and let $A$ be the set of $\alpha\in\mathcal A$
such that $x\in D_\alpha$. If $t_\alpha$ is a local equation of $D_\alpha$
in a neighborhood
$U$ of some $\R$-point of $D_\alpha$,
then there is a $\mathcal C^\infty$ function
$\phi_\alpha$ on $U$ such that
for $\mathbf x\in U\cap \ga^n(\R)$ we have
\[ h_\alpha(\mathbf x) = \log \abs{t_\alpha(\mathbf x)} +
\phi_\alpha(\mathbf x) . \]
It then follows from Propositions~\ref{lemm.vectorfields}
and~\ref{prop.logvanishing} that
for each $\alpha$,
\[ \frac{\partial}{\partial x_j} h_\alpha(\mathbf x) 
   = \frac12 \frac{\partial}{\partial x_j}
\log \abs{t_\alpha(\mathbf x)} + \frac{\partial}{\partial x_j} \phi_\alpha
(\mathbf x) \]
extends to a $\mathcal C^\infty$ function on $X(\R)$.

{}From the equality
\[ H_v(\mathbf x ;\mathbf s)^{-1}=\prod_{\alpha\in\mathcal A}
   \exp(-s_\alpha h_\alpha(\mathbf x)) \]
we deduce by induction the existence of an isobaric polynomial
$P_N\in\R[X^{(1)}_\alpha,
\dots, X^{(N)}_\alpha]$  of degree $N$ (with the convention that
each $X^{(p)}_\alpha$ has weight $p$) 
and such that
\[ (\partial/\partial x_j)^N H_v(\mathbf x ;\mathbf s)^{-1}
= H_v(\mathbf x ;\mathbf s)^{-1} P_N( s_\alpha \partial_j h_\alpha,
   s_\alpha\partial_j^2 h_\alpha, \dots, s_\alpha \partial_j^N h_\alpha).
\]
This implies that there exists a constant $C_v({\eps,N,j})$ such that
\[ \abs{(2i\pi a_j)^N} \abs{\hat H_v( \psi_{\mathbf a} ; \mathbf s )}
  \leq C_v({\eps,N,j}) (1+\norm{\mathbf s})^N \int_{\R^n}
      \abs{H_v(\mathbf x ;\mathbf s)}^{-1} \,\dx . \]
Choosing $j$ such that $\abs{a_j}$ is maximal  gives
$\abs{a_j}\geq \norm{\mathbf a}/\sqrt n$, hence an upper bound
\[ \abs{\hat H_v( \psi_{\mathbf a} ; \mathbf s )} 
  \leq C'_v(\eps,N) \frac{(1+\norm{\mathbf s})^N}{(1+\norm{\mathbf a})^N}
     \hat H_v(\Re(\mathbf  s);\psi_{0}),
\]
where $\psi_{0}$ is the trivial character and $C'_v(\eps, N)$
some positive constant. 
To conclude the proof it suffices to remark
that for any $\eps>0$,
$\hat H_v$ is bounded on the set $\Omega_{-1+\eps}$
(but the bound depends on $\eps$).

The case $F_v=\C$ is treated using a similar integration by parts.
(The exponent~$2$ on the numerator comes from the fact that
for a complex place $v$, $\norm{\cdot}_v$ is the square of a norm.)
\end{proof}

\section{Nonarchimedean computation at the trivial character}
\label{section.trivchar}

In this section we consider only $v\notin S$.
Let $\mathfrak m_v\subset\mathfrak o_v$
be the maximal ideal, $k_v=\mathfrak o_v/\mathfrak m_v$
and $q=\#k_v$.
Recall that we have fixed a good model $\mathcal X_{/U}$
over some $U\subset \Spec \mathfrak o_F$. 
To simplify notations  
we will write $\mathbf x=\mathbf x_v$, $\dx =\dx_v$ etc. 

The following formula is an analogue of Denef's formula
in~\cite[Thm~3.1]{denef87} for Igusa's local zeta function.
\begin{theo}\label{theo.trivcalc}
For all $v\notin S$ and all 
$\mathbf s \in \Omega_{-1}\subset \Pic(X)_{\C}$ 
we have
\[ \hat H_v(\psi_{0};\mathbf s) = q^{-\dim X}
         \sum_{A\subset\mathcal A} 
         \# D_{A}^\circ(k_v)
         \prod_{\alpha\in {A}} 
         \frac{q-1}{q^{1+s_{\alpha}-\rho_\alpha}-1}.
\]
\end{theo}

\begin{rema}
For $\mathbf s=\rho$ we get $\#  \mathcal X_{/U}(k_v)/q^{\dim X}$, the
expected local density at $v$.
\end{rema}

\begin{proof}
We split the integral along residue classes
mod $\mathfrak m_v$. Let $\tilde x\in \mathcal X(k_v)$
and $A=\{\alpha\in \mathcal A\,;\, \tilde x\in  D_{\alpha}\}$
so that $\tilde x\in  D_{A}^\circ$. 

We can introduce local ({\'e}tale) 
coordinates $x_{\alpha}$ ($\alpha\in A$)
and $y_{\beta}$ ($\beta\in B$, $\# A+\# B=\dim X$) around 
$\tilde x$
such  that locally, the divisor
$ D_{\alpha}$ 
is defined by the vanishing of $x_{\alpha}$.
Then, the local Tamagawa measure identifies with
the measure $\prod\mathrm dx_\alpha\times \prod \mathrm dy_\beta$
on $\mathfrak m_v^A\times\mathfrak m_v^B$. If $\mathrm d\mathbf x$
denotes the fixed measure on $\ga^n(F_v)$, one has the equality
of measures on $\ga^n(F_v)\cap \red^{-1}(\tilde x)$:
\[ \mathrm d\mathbf x = H_v(\mathbf x;\rho)\,\mathrm d\mu_v
   = \prod_{\alpha\in A} q^{\rho_\alpha v(x_\alpha)}\,
\prod\mathrm dx_\alpha\times \prod \mathrm dy_\beta.
\]
Consequently,
\begin{align*}
 \int_{\red^{-1}(\tilde x)} H_v(\mathbf x;\mathbf s)^{-1}\dx & = 
\int_{\mathfrak m_v^{A}\times\mathfrak m_v^{B}}
         q^{-\sum_{\alpha\in A} (s_{\alpha} -\rho_{\alpha}) v(x_{\alpha})}
         \prod_{\alpha\in A} \mathrm dx_{\alpha} 
         \prod_{\beta\in B}\mathrm dy_{\beta} \\
&= \frac1{q^{\# B}} \prod_{\alpha\in {A}} \int_{\mathfrak m_v} 
q^{-(s_{\alpha} -\rho_{\alpha}) v(x)}\,\mathrm dx 
= \frac1{q^{\dim X}} \prod_{\alpha\in A} 
          \frac{q-1}{q^{1+s_{\alpha}-\rho_{\alpha}}-1} 
\end{align*}
where the last equality follows from
\begin{equation}
\int_{\mathfrak m_v} q^{-sv(x)}\,\mathrm dx =
       \sum_{n=1}^\infty q^{-sn} 
\vol(\mathfrak m_v^n\setminus\mathfrak m_v^{n+1}) 
= \sum_{n=1}^\infty q^{-sn} q^{-n} \big(1-\frac1q\big) 
= \frac{1}q \frac{q-1}{q^{1+s}-1}.
\end{equation}
Summing these equalities for $\tilde x\in \mathcal X_{/U}(k_v)$
gives the desired formula.
\end{proof}

To estimate further $\hat H_v(\psi_0;\mathbf s)$ and the product
over all places,
we need to give an estimate for the number
of $k_v$-points in $D_{\alpha}$ which is uniform in  $v$.

\begin{lemm}\label{lemm.uniform}
There exists a constant $C(X)$ such that for all $v\notin S$ 
and all $A\subset \mathcal A$ 
we have the estimates:
\begin{itemize}
\item if $\#A =1$,
$
 \abs{\#  D_{A}(k_v) - q_v^{\dim X -1}} \le C(X)q_v^{\dim X - 3/2}
$;
\item 
if $\#A \ge 2$,
$ \# D_{A}(k_v)\le C(X)q_v^{\dim X - \#A}$.
\end{itemize}
\end{lemm}
\begin{proof}
We use the fact 
that for any projective
variety $Y$ of dimension $\dim Y$ and degree $\deg Y $ the number of 
$k_v$-points is estimated as
\[
\#Y(k_v)\le (\deg Y) \#\P^{\dim Y}(k_v),
\]
(see, e.g., Lemma 3.9 in \cite{chambert-loir-t2000b}). 
Since $X$ is projective, all $D_{A}$  
can be realized as subvarieties 
in some projective space. This proves the second part. 
To prove the first part, we apply Lang-Weil's estimate~\cite{lang-w54} 
to the geometrically irreducible 
smooth $U$-scheme $ D_{\alpha}$. 
\end{proof}

\begin{prop}\label{prop.bbb}
For all  $\eps>0$  there exists a constant $C(\eps)$
such that for any $\mathbf s\in\Omega_{-\frac12+\eps}$ and 
any finite place $v\not\in S$,
\[ \abs{ \hat H_v(\psi_{0};\mathbf s)
   \prod_{\alpha\in \mathcal A}(1-q^{-(1+s_{\alpha}-\rho_{\alpha})})
   - 1 } 
   \leq C(\eps) q^{-1-\eps}. \]
\end{prop}
\begin{proof}
Using the uniform estimates from Lemma~\ref{lemm.uniform} 
we see that in the formula for $\hat H_v$, each term with $\# A\geq 2$
is $O(q^{-(\frac12+\eps)\#A})=O(q^{-1-2\eps})$, with uniform 
constants in $O$.
Turning to the remaining terms, we get 
\begin{align*}
 & \hphantom{{}={}} 1+\sum_{\alpha\in \mathcal A}
     (\frac 1q +O(1/q^2))
     \frac {q-1}{q^{1+s_{\alpha}-\rho_{\alpha}}-1} \\
& = 1+ \sum_{\alpha\in \mathcal A}
             q^{-(1+s_{\alpha}-\rho_{\alpha})} (1-\frac 1q)
             (1-q^{-(1+s_{\alpha}-\rho_{\alpha})})^{-1}
         + O(q^{-3/2}) \\
&= 1+ \sum_{\alpha\in \mathcal A}
       \frac{q^{-(1+s_{\alpha}-\rho_{\alpha})}}
            {1-q^{-(1+s_{\alpha}-\rho_{\alpha})}}  + O(q^{-3/2}) \\
& = \prod_{\alpha\in \mathcal A}
           (1-q^{-(1+s_{\alpha}-\rho_{\alpha})})^{-1}
           (1+O(q^{-1-2\eps})) + O(q^{-3/2}) .
\end{align*}
Finally, we have the desired estimate.
\end{proof}

For $X$ as above and $\mathbf s=(s_{\alpha})\in \Pic(X)_{\C}$,
the (multi-variable) Artin $L$-function is given by
\[
L(\Pic(X);\mathbf s) = \prod_{\alpha\in\mathcal A}\zeta_F(s_\alpha)=
  \prod_{\text{$v$ finite}}
  \prod_{\alpha\in\mathcal A} (1-q^{-s_{\alpha}})^{-1} 
.
\]
From standard properties of Dedekind zeta functions,
we conclude that $L(\Pic(X);\mathbf s)$ admits a meromorphic continuation
and that it has polynomial growth in vertical strips. 
For $\mathbf s=(s,\dots,s)$ we get the usual Artin $L$-function
of $\Pic(X)$---here a power of the Dedekind zeta function---which
has been used in the regularization of the Tamagawa measure
in Prop.~\ref{prop.tamagawa}.
 
\begin{coro}\label{coro.trivtamagawa}
For all $\eps >0$ there 
exists a holomorphic bounded function $\phi(0;\cdot)$
on $\Omega_{-1/2 + \eps}$ such
that for any $\mathbf s\in \Omega_0$ one has
\[ \hat H(\psi_{0};\mathbf s) = 
\phi(0;\mathbf s) \prod_{\alpha\in\mathcal A}(s_\alpha-\rho_\alpha)^{-1}.
\] 
Moreover, there exist constants $N>0$ and 
$C(\eps)$ such that for any $\mathbf s$ in $\Omega_{-1/2+\eps}$,
\[ \abs{\phi(0;\mathbf s)} \leq C(\eps) (1+\norm{\Im(\mathbf s)})^{N}.
\]
\end{coro}

\begin{proof}
For any place~$v$ of~$F$, let
\[ f_v (\mathbf s)
= \hat H_v(\psi_0;\mathbf s) L_v(\Pic(X);\mathbf s-\rho+1)^{-1}. \]
The previous proposition shows that in $\Omega_{-1/2+ \eps}$,
the Euler product $\prod_v f_v$
converges absoluetely to a holomorphic and bounded function~$f$.
For any $\mathbf s\in\Omega_0$, one has
\[ \hat H(\psi_{0};\mathbf s)
 = \prod_v \hat H_v(\psi_0;\mathbf s)
 = \prod_v f_v(\mathbf s) \prod_v L_v(\Pic(X);\mathbf s-\rho+1)
 = f(\mathbf s) L(\Pic(X);\mathbf s).
\]
It suffices to define
\[ \phi(0;\mathbf s) = f(\mathbf s) L(\Pic(X);\mathbf s-\rho+1)
 \prod_{\alpha\in\mathcal A}(s_\alpha-\rho_\alpha).
\]
The polynomial growth of $\phi(0;\cdot)$ in vertical strips
follows from the boundedness of~$f$ in $\Omega_{-1/2+\eps}$
and from  the fact that Dedekind zeta functions have 
a polynomial growth in such vertical strips.
\end{proof}

\section{Other characters} 
\label{section.otherchars}

Our aim here is to compute as explicitely as possible
the Fourier transforms of local heights with
character $\psi_{\mathbf a}$. In general, this will be only possible
up to some error term.

The calculations in the preceding Section allow us to strengthen
Proposition~\ref{prop.estim}. For $\mathbf a\in\ga^n(F)$, 
we denote by $\norm{\mathbf a}_\infty$ a norm of $\mathbf a$
in the real  vector space
defined by extension of scalars
via the diagonal embedding $F\hookrightarrow F\otimes_{\Q}\R
= \prod_{v\mid\infty}F_v$.

\begin{prop}
For any $\eps>0$, there exist an integer $\kappa\geq 0$
and for any $N\geq 0$,
a constant  $C(\eps,N)$ such that for any
$\mathbf s\in \Omega_{-1/2 +\eps}$ and all
${\mathbf a}\in \mathfrak d_X$  we have 
\[ \prod_{v\in S(\mathbf a)} \abs{\hat H_v( \psi_{\mathbf a} ; \mathbf s )}
   \le C(\eps,N) (1+\norm{\mathbf s})^{N}
                 (1+\norm{\mathbf a}_\infty)^{\kappa-N}. \]
\end{prop}
\begin{proof}
For $v\in S$ the local integrals converge absolutely in 
the domain under consideration and are bounded as in 
Proposition~\ref{prop.estim}.
For $v\not\in S$, 
we have shown that
there exists a constant $c$ 
such that for all $\mathbf a \in \mathfrak d_X$
and all $\mathbf s \in \Omega_{-1/2+\eps}$ 
one has the estimate
\[
\abs{\hat H_v( \psi_{\mathbf a} ; \mathbf s )-1}\le \frac{c}{q_v^{\eps}}.
\] 
This implies that 
$\abs{\hat H_v(\psi_ {\mathbf a};\mathbf s)}$ is bounded
independently of $\mathbf a$, $v\in S(\mathbf a)\setminus S$ and 
$\mathbf s\in \Omega_{-1/2+\eps}$.
For any nonzero $a\in\mathfrak o_F$,
there is a trivial estimate 
\[ \sum_{\mathfrak p\supset (a)} 1 
\ll \sum_{\mathfrak p\supset (a)} \log \mathcal N(\mathfrak p)
\ll \log \mathcal N(a), \] 
which implies that
for some constant $\kappa$,
\[ \prod_{v\in S(\mathbf a)\setminus S}
	\abs{\hat H_v(\psi_ {\mathbf a};\mathbf s)}
 \ll \prod_{v|\infty} (1+\norm{\mathbf a}_v)^\kappa.\]
Using Proposition~\ref{prop.estim}, we have for all $N> 0$,
\begin{align*}
 \prod_{v\in S(\mathbf a)} \abs{\hat H_v( \psi_{\mathbf a} ; \mathbf s )}
 & \ll C(S,\eps,N) (1+\norm{s})^{N[F:\Q]} \prod_{v|\infty} (1+\norm{\mathbf a}_v)^{\kappa-N} \\
 & \leq C'(S,\eps,N) (1+\norm{s})^{N[F:\Q]} (1+\norm{\mathbf a}_\infty)^{(\kappa-N)[F:\Q]}
.\end{align*}
Replacing $N$ by $N[F:\Q]$ and $\kappa$
by $\kappa[F:\Q]$ concludes the proof of the proposition.
\end{proof}

For the explicit calculation at good places,
let us recall some notations: 
$\mathfrak m_v\subset \mathfrak o_v$ is the maximal ideal, 
$\pi=\pi_v$ a uniformizing element, 
$k_v$ the residue field, $q=q_v=\#k_v$, $\psi =\psi_v=\psi_{\mathbf a,v}$, 
$\mathbf x=\mathbf x_v$,
$\mathbf a =\mathbf a_v$, $\dx =\dx_v$.     
As in Section~\ref{section.geometry}, let $f=f_{\mathbf a}$ be a linear
form on $\ga^n$ and write
$\div(f)=E-\sum_\alpha d_\alpha D_\alpha$.
We have also defined
$\mathcal A_0(\mathbf a)=\{\alpha\,;\, d_\alpha(f_{\mathbf a})=0\}$
and
$\mathcal A_1(\mathbf a)=\{\alpha\,;\, d_\alpha(f_{\mathbf a})=1\}$.

\begin{prop}\label{prop.estimatewithcharacters}
There exists a constant $C(\eps)$ independent of 
$\mathbf a\in \mathfrak d_X$
such that for any $v\not\in S(\mathbf a)$,
\[  \abs{ \hat H_v( \psi_{\mathbf a} ; \mathbf s )
 \prod_{\alpha\in \mathcal A_0(\mathbf a)} 
(1-q^{-(1+s_\alpha-\rho_\alpha)}) - 1 }
     \leq C(\eps) q^{-1-\eps}. \]
\end{prop}

Similarly to 
the proof of Proposition~\ref{prop.bbb} in the preceding Section,
this proposition is proved by computing the integral on residue classes.

Let $\tilde x\in X(k_v)$ and $A=\{\alpha\,;\, 
\tilde x\in D_{\alpha}\}$.
We now consider three cases:

\emph{Case 1. $A=\emptyset$.} ---
The sum of all these contribution
is equal to the integral over $\ga^n(\mathfrak o_v)$:
\[
\int_{\ga(\mathfrak o_v)} H_v(\mathbf x;\mathbf s)^{-1}
     \psi_{\mathbf a}(\mathbf x)\, \mathrm d\mathbf x 
= \int_{\mathfrak o_v^n} \psi( \langle \mathbf a,\mathbf x\rangle)
       \, \mathrm d\mathbf x  
= 1.
\]

\emph{Case 2. $A=\{\alpha\}$ and $\tilde x\not\in E$.} ---
We can introduce analytic coordinates $x_{\alpha}$ and $y_{\beta}$
around $\tilde x$ such that locally
$f(x)=u x_{\alpha}^{-d_{\alpha}}$ with $u\in\mathfrak o_v^*$.
Then, we compute the integral of  
$H_v(\mathbf x ;\mathbf s)^{-1}\psi_{\mathbf a}(\mathbf x)$ as 
\begin{align*}
\int_{\red^{-1}(\tilde x)} &= 
   \int_{\mathfrak m_v\times\mathfrak m_v^{n-1}}
        q^{-(s_\alpha-\rho_\alpha)v(x_\alpha)} 
        \psi(u x_{\alpha}^{-d_{\alpha}})
         \, \mathrm dx_{\alpha}
         \mathrm d\mathbf y \\
&= \frac{1}{q^{n-1}}
   \sum_{n_{\alpha}\geq 1}
            q^{- (1+s_{\alpha}-\rho_{\alpha})n_{\alpha}}
            \int_{\mathfrak o_v^*}
              \psi(u \pi^{- n_{\alpha} d_{\alpha}} u_{\alpha}^{-d_{\alpha}})\,
               \mathrm du_{\alpha} .
\end{align*}

\begin{lemm}\label{lemm.cancel}
For all integers $d\geq 0$ and $n\geq 1$ and all $u\in\mathfrak o_v^*$,
\[ \int_{\mathfrak o_v^*} \psi(u\pi^{-nd}t^d)\,\mathrm dt
      = \begin{cases}
       1-1/q & \text{if $d=0$;} \\
        -1/q & \text{if $n=d=1$;} \\
        0 &\text{else.}
        \end{cases}
\]
\end{lemm}
\begin{proof}
If $d=1$, the computation runs as follows
\begin{align*}
\int_{\mathfrak o_v^*} \psi(u\pi^{-n}t)\,\mathrm dt 
&= \int_{\mathfrak o_v} \psi(u\pi^{-n} t)\,\mathrm dt
   - \frac1q \int_{\mathfrak o_v} \psi(u\pi^{-n+1}t)\,\mathrm dt \\
&= \begin{cases}
    0 & \text{if $n\geq 2$;} \\
    -1/q & \text{if $n=1$.} 
  \end{cases}
\end{align*}
For $d\geq 2$,
let $r=v_p(d)$. Since we assumed $F_v$ to be unramified
over~$\Q_p$, $r=v_{\pi}(d)$.
We will integrate over disks $D(\xi,\pi^e)\subset \mathfrak o_v^*$
for $e\geq 1$ suitably chosen. Indeed, if $v\in\mathfrak o_v$
and $t=\xi +\pi^e v$,
\[ t^d = \xi^d+ d\pi^e \xi^{d-1} v \pmod {\pi^{2e} } \]
hence, if $a$ is chosen such that
\[ e-nd+r <0 \quad\text{and}\quad 2e-nd \geq 0, \]
\[ \int_{D(\xi,\pi^e)} \psi(u\pi^{-nd}t^d)\, \mathrm dt
  = q_v^{-e} \psi(u\pi^{-nd}\xi^d)
             \int_{\mathfrak o_v} \psi(d\pi^{e-nd} u\xi^{d-1}\, v) \,
\mathrm dv = 0.
\]
We can find such an $e$ if and only if $2(nd-r-1)\geq nd$, \emph{i.e.}
$nd \geq 2r+2$.
If $r=0$, this is true since $d\geq 2$. If $r\geq 1$, one has
$nd \geq p^r \geq 2r+2$ since we assumed $p\geq 5$.
\end{proof}

This lemma implies  the following trichotomy:
\begin{align*}
 \int_{\red^{-1}(\tilde x)} 
H(\mathbf x;\mathbf s)^{-1} \psi_{\mathbf a}(\mathbf x) \, \mathrm d\mathbf x
& = \frac{q-1}{q^n} \frac{1}{q^{1+s_\alpha-\rho_\alpha}-1}  & \text{if
$d_\alpha=0$;} \\
& = -\frac{1}{q^n} q^{-(1+s_\alpha-\rho_\alpha)} & \text{if $d_\alpha=1$;} \\
& = 0 & \text{if $d_\alpha\geq 2$.} \\
\end{align*}

\emph{Case 3. $\# A\geq 2$ or $\# A=1$ and $\tilde x\in E$.} ---
Under some transversality assumption, we could compute explicitely
the integral as before. We shall however content ourselves with
the estimate obtained by replacing $\psi$ by 1 in the integral.

The total contribution of these points
will therefore be smaller than
\begin{equation} 
\label{eq.ET}
 q^{-\dim X} \sum_{\# A\geq 2} \# D_A^\circ(k_v) \prod_{\alpha\in A}
    \frac{q-1}{q^{1+s_\alpha-\rho_\alpha}-1}
  + q^{-\dim X} \sum_{A=\{\alpha\}} \#(D_\alpha\cap E)(k_v)
         \frac{q-1}{q^{1+s_\alpha-\rho_\alpha}-1}. 
\end{equation}

Finally, the Fourier transform 
is estimated as follows:
\[ 
 \hat H_v(\psi_{\mathbf a};\mathbf s)
= 1+ q^{-n} \sum_{\alpha\in \mathcal A_0(\mathbf a)}
     \# D_\alpha^\circ (k_v) \frac{q-1}{q^{1+s_\alpha-\rho_\alpha}-1} 
{} - q^{-n} \sum_{\alpha\in \mathcal A_1(\mathbf a)}
   \# D_\alpha^\circ(k_v) \frac{1}{q^{1+s_\alpha-\rho_\alpha}}
 +  \mathit{ET} 
\] 
with an ``error term''~$\mathit{ET}$
smaller than the expression in~\eqref{eq.ET}.
It is now a simple matter to rewrite this estimate as in the statement
of~\ref{prop.estimatewithcharacters}.\qed

\medskip

We deduce from these estimates that 
$\hat H(\psi_{\mathbf a};\mathbf s)$ has a meromorphic
continuation:
\begin{coro}\label{coro.otherchars}
For any $\eps >0$ and $\mathbf a\in \mathfrak d_X \setminus \{ 0\}$
there exists a holomorphic bounded function $\phi(\mathbf a;\cdot)$
on $\Omega_{-1/2+\eps}$ such
that for any $\mathbf s\in \Omega_{0}$ 
\[ \hat H(\psi_{\mathbf a};\mathbf s)=\prod_v \hat H_v( \psi_{\mathbf a} ; \mathbf s ) = 
\phi( \psi_{\mathbf a} ; \mathbf s )\prod_{\alpha\in \mathcal A_0(\mathbf a)}
     \zeta_F(1+s_\alpha-\rho_\alpha). \]
Moreover, for any $N>0$ there exist  
constants $N'>0$ and $C(\eps, N)$ such that
for any $\mathbf s\in\Omega_{-1/2+\eps}$,
one has the estimate
\[ \abs{\phi(\psi_{\mathbf a};\mathbf s)} \leq C(\eps,N) 
     (1+\norm{\Im(\mathbf s )})^{N'}
     (1+\norm{\mathbf a}_\infty)^{-N}. \]
\end{coro}

\section{The nonsplit case}\label{section.nonsplit}

In this section we extend the previous calculations of the Fourier
transform to the nonsplit case, \emph{i.e.} when the geometric irreducible
components of $X\setminus\ga^n$ are no longer assumed to be
defined over~$F$.

\begin{lemm}
Let $x\in X(F)$ and $A=\{\alpha\in\mathcal A\,;\, x\in D_\alpha\}$.
Fix for any orbit $\bar\alpha \in A/\Gamma_F$ some element $\alpha$
and let $F_\alpha$ be the field of definition of $D_\alpha$.
Then there exist an open neighborhood~$U$ 
of~$x$, {\'e}tale coordinates around~$x$ over $\bar F$,
$x_\alpha$ ($\alpha\in A$) and $y_\beta$ such that 
$x_\alpha$ is a local equation of $D_\alpha$ and such
that the induced morphism $U_{\bar F}\ra \A^n_{/\bar F}$
descends to an {\'e}tale morphism
\[ U \ra \prod_{\bar\alpha\in A/\Gamma_F} \Res_{F_\alpha/F} \A^1
          \times \A^{n-a} \qquad (a=\# A). \]

An analogous result holds over the local fields~$F_v$ and also
on $\mathfrak o_v$, $v$ being any finite place of $F$  such that
$v\not\in S$.
\end{lemm}
(We have denoted by $\Res$ the functor of restriction of scalars
\emph{\`a la} Weil.)
\begin{proof}
Chose some element $\alpha$ in each orbit $\bar\alpha$
and fix
a local equation $x_\alpha$ for the corresponding $D_\alpha$
which is defined over the number field $F_\alpha$.
Then if $\alpha'=g\alpha$ (for some $g\in \Gamma_F$)
is another element in the orbit of~$\alpha$,
set $x_{\alpha'}=g\cdot x_\alpha$. This is well defined since
we assumed that the equation $x_\alpha$ is invariant under
$\Gamma_{F_\alpha}$.

Finally, add $F$-rational local  {\'e}tale coordinates corresponding to
a basis of the subspace in $\Omega^1_{X,x}$ which is complementary to
the span of $dx_\alpha$ for $\alpha\in A$.
\end{proof}

Let $v$ be a place of~$F$.
The above lemma allows us to identify a neighborhood of~$x$
in $X(F_v)$
(for the analytic topology) with a neighborhood of~$0$ in the product
$\prod_{\alpha\in  A/\Gamma_v} F_{v,\alpha} \times F_v^{n-a}$,
the local heights $\prod_{\alpha\in\alpha} H_{\alpha,v}(\xi)$
being replaced by 
\[ \mathcal N_{F_{v,\alpha}/F_v}(\xi) \times h_{\alpha,v}(\xi) \]
where $h_{\alpha,v}$ is a smooth function.

Similarly, for good places~$v$ we identify $\red^{-1}(\tilde x)$
with $\prod {\mathfrak m_{v,\alpha}} \times \mathfrak m^{n-a}$
and the functions $h_{\alpha,v}$ are equal to~$1$.

The assertions of~Section~\ref{section.estimates}  
still hold in this more general case
and the proofs require only minor modifications.
However the calculations of~Sections~\ref{section.trivchar}
and~\ref{section.otherchars} have to be redone.

\begin{theo}[cf.~Thm.~\ref{theo.trivcalc}]
One has
\[ \hat H_v(\psi_0;\mathbf s)
      = q_v^{-\dim X}
        \sum_{A\subset \mathcal A/\Gamma_v}
        \# D_A^\circ(k_v)
        \prod_{\alpha \in A/\Gamma_v}
           \frac{q_v^{f_\alpha}-1}{q_v^{f_\alpha(1+s_\alpha-\rho_\alpha)}-1}            
\]
where $f_\alpha$ is degree of $F_{v,\alpha}$ over $F_v$.
\end{theo}

\begin{coro}[cf.~Prop.~\ref{prop.bbb}]
One has
\[ \hat H_v(\psi_0;\mathbf s) =
      \prod_{\alpha\in \mathcal A/\Gamma_v}
         \big( 1-q_v^{-f_\alpha(1+s_\alpha-\rho_\alpha)}\big)^{-1}
        \big( 1+ O(q_v^{-1-\eps}) \big). 
\]
\end{coro}

In the general case, we introduce
the multi-variable Artin $L$-function of $\Pic(\bar X)$ 
as
\[ L(\Pic(\bar X);\mathbf s) =  \prod_{\alpha\in\mathcal A/\Gamma_F}
          \zeta_{F_\alpha}(s_\alpha). \]
Its restriction to the line $(s,\dots,s)$ is the usual Artin $L$-function
of $\Pic(\bar X)$. It has a pole of order $\#(\mathcal A/\Gamma_F)
= \rk(\Pic X)$ at $s=1$.

\begin{coro}[cf.~Cor.~\ref{coro.trivtamagawa}]
For all~$\eps>0$, there exists a holomorphic bounded function $\phi$
on $\Omega_{-1/2+\eps}$ such that for any $\mathbf s\in\Omega_0$ one has
\[ \hat H(\psi_0;\mathbf s) = \phi(\mathbf s) L(\Pic(\bar X);\mathbf
s-\rho+1) \]
and the Tamagawa measure of~$X(\A_F)$ is equal to
\[ \tau(K_X) = \phi(\rho)/ L^*(\Pic(\bar X);\mathbf 1). \]
\end{coro}

At nontrivial characters, the calculations are modified analogously
and we have
\[ \hat H(\psi_{\mathbf a};\mathbf s)
     = \phi(\psi_{\mathbf a};\mathbf s) \prod_{\alpha\in \mathcal A_0(\mathbf a)/\Gamma_F}
         \zeta_{F_\alpha}(s_\alpha-\rho_\alpha+1) 
\]
for some holomorphic function $\phi(\mathbf a;\cdot)$ 
as in~Corollary~\ref{coro.otherchars}.

\def\noop#1{\ignorespaces}
\nocite{chambert-loir-t2000,chambert-loir-t2000b,tate67b}
\nocite{peyre95,batyrev-m90,batyrev-t98b,batyrev-t98}


\providecommand{\noopsort}[1]{}
\providecommand{\bysame}{\leavevmode ---\ }
\providecommand{\og}{``}
\providecommand{\fg}{''}
\providecommand{\smfandname}{et}
\providecommand{\smfedsname}{\'eds.}
\providecommand{\smfedname}{\'ed.}
\providecommand{\smfmastersthesisname}{M\'emoire}
\providecommand{\smfphdthesisname}{Th\`ese}

\end{document}